\documentclass{amsart}

\usepackage{amsfonts,mathrsfs}
\usepackage{amsmath,amssymb,amsthm, amscd, bm}
\usepackage{tikz}
\usetikzlibrary{arrows,calc,matrix,topaths,positioning,scopes,shapes,decorations,decorations.markings} 
\usepackage{dsfont}
\usepackage{epsfig}
\usepackage{tikz-cd}
\usepackage{hyperref}
\usepackage{fullpage}
\usepackage[foot]{amsaddr}
\usepackage{adjustbox}
\usepackage{enumitem}

\usepackage[utf8x]{inputenc}

\makeatletter

\renewcommand{\email}[2][]{%
  \ifx\emails\@empty\relax\else{\g@addto@macro\emails{,\space}}\fi%
  \@ifnotempty{#1}{\g@addto@macro\emails{\textrm{(#1)}\space}}%
  \g@addto@macro\emails{#2}%
}
\makeatother

\author{Julien Korinman}
\address{ Institut Montpelli\'erain Alexander Grothendieck - UMR 5149 Universit\'e de Montpellier. Place Eug\'ene Bataillon, 34090 Montpellier France}
\email{julien.korinman@gmail.com}
\keywords{Stated skein algebras, mapping class groups}

\def\restriction#1#2{\mathchoice
              {\setbox1\hbox{${\displaystyle #1}_{\scriptstyle #2}$}
              \restrictionaux{#1}{#2}}
              {\setbox1\hbox{${\textstyle #1}_{\scriptstyle #2}$}
              \restrictionaux{#1}{#2}}
              {\setbox1\hbox{${\scriptstyle #1}_{\scriptscriptstyle #2}$}
              \restrictionaux{#1}{#2}}
              {\setbox1\hbox{${\scriptscriptstyle #1}_{\scriptscriptstyle #2}$}
              \restrictionaux{#1}{#2}}}
\def\restrictionaux#1#2{{#1\,\smash{\vrule height .8\ht1 depth .85\dp1}}_{\,#2}}

\newcommand{\quotient}[2]{{\raisebox{.2em}{$#1$}\left/\raisebox{-.2em}{$#2$}\right.}}

\newcommand{\sslash}{\mathbin{/\mkern-6mu/}}
\newcommand{\Hom}{\operatorname{Hom}}
\newcommand{\tr}{\operatorname{tr}}

\newcommand{\qtr}{\operatorname{qtr}}
\newcommand{\SL}{\operatorname{SL}}
\newcommand{\id}{id}
\newcommand{\Span}{\operatorname{Span}}
\newcommand{\End}{\operatorname{End}}
\newcommand{\GL}{\operatorname{GL}}
\newcommand{\PGL}{\operatorname{PGL}}
\newcommand{\Vect}{\operatorname{Vect}}
\newcommand{\Specm}{\operatorname{MaxSpec}}

\newcommand{\Mod}{\operatorname{Mod}}

\newcommand{\Mat}{\operatorname{Mat}}

\newcommand{\Comod}{\operatorname{Comod}}

\newcommand{\MS}{\operatorname{MS}}
\newcommand{\Alg}{\operatorname{Alg}}
\newcommand{\Irrep}{\operatorname{Irrep}}

\newcommand{\Proj}{Proj}

\newcommand{\Cob}{Cob}
\newcommand{\BT}{BT}

\newcommand{\heightexch}[3]{
	\begin{tikzpicture}[baseline=-0.4ex,scale=0.5, >=stealth]
	\draw [fill=gray!60,gray!45] (-.7,-.75)  rectangle (.4,.75)   ;
	\draw[#1] (0.4,-0.75) to (.4,.75);
	\draw[line width=1.2] (0.4,-0.3) to (-.7,-.3);
	\draw[line width=1.2] (0.4,0.3) to (-.7,.3);
	\draw (0.65,0.3) node {\scriptsize{$#2$}}; 
	\draw (0.65,-0.3) node {\scriptsize{$#3$}}; 
	\end{tikzpicture}
}
\newcommand{\heightcurve}{
\begin{tikzpicture}[baseline=-0.4ex,scale=0.5]
\draw [fill=gray!20,gray!45] (-.7,-.75)  rectangle (.4,.75)   ;
\draw[-] (0.4,-0.75) to (.4,.75);
\draw[line width=1.2] (-.7,-0.3) to (-.4,-.3);
\draw[line width=1.2] (-.7,0.3) to (-.4,.3);
\draw[line width=1.15] (-.4,0) ++(-90:.3) arc (-90:90:.3);
\end{tikzpicture}
}

\begin{document}

\theoremstyle{plain}
\newtheorem{theorem}{Theorem}[section]
\newtheorem{proposition}[theorem]{Proposition}
\newtheorem{corollary}[theorem]{Corollary}
\newtheorem{lemma}[theorem]{Lemma}
\theoremstyle{definition}
\newtheorem{notations}[theorem]{Notations}
\newtheorem{convention}[theorem]{Convention}
\newtheorem{problem}[theorem]{Problem}
\newtheorem{definition}[theorem]{Definition}
\theoremstyle{remark}
\newtheorem{remark}[theorem]{Remark}
\newtheorem{conjecture}[theorem]{Conjecture}
\newtheorem{example}[theorem]{Example}
\newtheorem{strategy}[theorem]{Strategy}
\newtheorem{question}[theorem]{Question}

\title[A survey on mapping class group representations derived from skein theory]{A survey on mapping class group representations derived from skein theory}
%
%
%\author{Julien Korinman}
%

\date{}
\maketitle

%%%%%%%%%%%%%%%%%%%%%%%%%%%%%%%%%%%%%%%%%%%%%%%%%%%%%%%%%%%%%%%%%%%%%%%%%%%%%%%%

\begin{abstract} 
We survey various constructions of  finite dimensional projective representations of mapping class groups derived from stated skein algebras.
\end{abstract}

%\tableofcontents

%%%%%%%%%%%%%%%%%%%%%%%%%%%%%%%%%%%%%%%%%%%%%%%%%%%%%%%%%%%%%%%%%%%%%%%%%%%%%%%

\section{Overview}
Stated skein algebras are 	associative algebras $\mathcal{S}_A(\mathbf{\Sigma})$ associated to a complex number $A^{1/2}$ and a marked surface $\mathbf{\Sigma}$. Informally, they are defined as the quotient $\mathcal{S}_A(\mathbf{\Sigma})= \quotient{\Span_{\mathbb{C}}\left( (T,s) \subset \Sigma \times [0,1]\right)}{(\mbox{Skein relations})}$ of the linear span of  "stated tangles" $(T,s)$ embedded in the thickened surface $\Sigma\times [0,1]$ by the subspace spanned by the so-called "skein relations". The mapping class group $\Mod(\mathbf{\Sigma})$ naturally acts on $\mathcal{S}_A(\mathbf{\Sigma})$ by acting on the tangles. An irreducible representation $r: \mathcal{S}_A(\mathbf{\Sigma}) \to \End(V)$ will be called \textbf{fixed by} $\Mod({\Sigma})$ if for every $\phi\in \Mod(\mathbf{\Sigma})$, the representation $\phi\cdot r : \mathcal{S}_A(\mathbf{\Sigma}) \to \End(V), x \mapsto r( \phi^{-1}(x))$ is isomorphic to $r$. In this case, there exists an operator $L_r(\phi) \in \End(V)$, unique up to multiplication by a non zero scalar, such that 
$$ r(\phi^{-1}(x)) = L_r(\phi) r(x) L_r(\phi)^{-1}, \quad \mbox{ for all } x \in \mathcal{S}_A(\mathbf{\Sigma}).$$
The assignation $\phi \mapsto L_r(\phi)$ defines a projective representation $L_r: \Mod({\Sigma}) \to \PGL(V)$; we therefore have a procedure $r\mapsto L_r$ to derive mapping class group representations from certain skein representations. Similarly, if $G\subset \Mod(\mathbf{\Sigma})$ is a subgroup, such as the Torelli group, and  $r: \mathcal{S}_A(\mathbf{\Sigma}) \to \End(V)$  is an irreducible representation which is fixed by $G$, we obtain a projective representation $L_r : G \to \PGL(V)$. We will present a lot of families of representations which arise that way, with an emphasize on the following four families: 
\begin{itemize}
\item the Witten-Reshetikhin-Turaev representations  from \cite{RT}; 
\item the Blanchet-Costantino-Geer-Patureau Mirand representations \cite{BCGPTQFT}: these are representations of the Torelli group which were conjectured to be faithful; 
\item the representations of the braid groups which underlie quantum invariants of knots and links \cite{KojuQGroupsBraidings} such as: the colored Jones polynomials, the Kashaev invariant, the ADO invariants, the Geer-Patureau Mirand invariants and the Blanchet-Geer-Patureau Mirand-Reshetikhin invariants; 
\item the Lyubashenko representations from \cite{Lyubashenko_InvariantsMCGRep}.
\end{itemize}
Besides these four families, we will also present families of representations of $\Mod(\Sigma)$ which have never been considered in literature but might deserve the same attention.

\section{Stated skein algebras in a nutshell}

\underline{\textbf{Marked surfaces:}} A \textbf{marked surface} $\mathbf{\Sigma}=(\Sigma, \mathcal{A})$ is a compact oriented surface $\Sigma$ with a finite set $\mathcal{A}$ of closed pairwise disjoint arcs in its boundary named \textbf{boundary edges}. For $a,b$ two  
boundary edges, we denote by $\mathbf{\Sigma}_{a\# b}$ the marked surface obtained by gluing $a$ to $b$. An \textbf{embedding} $f: \mathbf{\Sigma} \to \mathbf{\Sigma}'$ is an oriented embedding $f:\Sigma \hookrightarrow \Sigma'$  of the underlying marked surfaces which restrict to an embedding $f: \mathcal{A} \hookrightarrow \mathcal{A}'$. When several boundary edges $a_1, \ldots, a_n \in \mathcal{A}$ are embedded to the same boundary edge $a' \in \mathcal{A}'$, we specify a total ordering of $\{a_1, \ldots, a_n\}$ as part of the definition of $f$. 
 $f$ is a \textbf{strong embedding} if for every $a' \in \mathcal{A}'$ there exists at most one $a \in \mathcal{A}$ which is embedded into $a'$. Marked surfaces with embeddings and strong embeddings form categories $\MS^{str} \subset \MS$. 
\par \underline{\textbf{Punctures:}} A \textbf{puncture} of $\mathbf{\Sigma}$ is a connected component of $\partial \Sigma \setminus \mathcal{A}$: it is called an \textbf{inner puncture} if it is a circle (an unmarked boundary component) and a \textbf{boundary puncture} if it is an open arc. When we draw pictures, we will depict marked surface by drawing the underlying surface and drawing a puncture $\bullet$ at the location of each punctures.

\par \underline{\textbf{Stated tangles:}} A \textbf{tangle} in $\mathbf{\Sigma}$  is a framed submanifold $T\subset \Sigma\times (0,1)$ such that $\partial T \subset \mathcal{A} \times [0,1]$. We impose the framing of each point of $\partial T$ to be in the $(0,1)$ direction towards $1$ and for each $a \in \mathcal{A}$, no two points of $\partial_a T:= T\cap a\times (0,1)$ can have the same height (projection in the $[0,1]$ factor). A \textbf{state} is a map $s: \partial T \to \{-, +\}$. A \textbf{stated tangle} is a pair $(T,s)$.

\begin{definition}(Stated skein algebras)\label{def_sskein}
Let $k$ be a unital commutative ring, $A^{1/2} \in k^*$ an invertible element and $\mathbf{\Sigma}$ a marked surface. 
The \textbf{stated skein algebra}  $\mathcal{S}_A(\mathbf{\Sigma})$ is the  quotient of the $k$ module freely generated by isotopy classes of stated tangles in $\mathbf{\Sigma}$ modulo the following skein relations: 
$$
\begin{tikzpicture}[baseline=-0.4ex,scale=0.5,>=stealth]	
\draw [fill=gray!45,gray!45] (-.6,-.6)  rectangle (.6,.6)   ;
\draw[line width=1.2,-] (-0.4,-0.52) -- (.4,.53);
\draw[line width=1.2,-] (0.4,-0.52) -- (0.1,-0.12);
\draw[line width=1.2,-] (-0.1,0.12) -- (-.4,.53);
\end{tikzpicture}
=A
\begin{tikzpicture}[baseline=-0.4ex,scale=0.5,>=stealth] 
\draw [fill=gray!45,gray!45] (-.6,-.6)  rectangle (.6,.6)   ;
\draw[line width=1.2] (-0.4,-0.52) ..controls +(.3,.5).. (-.4,.53);
\draw[line width=1.2] (0.4,-0.52) ..controls +(-.3,.5).. (.4,.53);
\end{tikzpicture}
+A^{-1}
\begin{tikzpicture}[baseline=-0.4ex,scale=0.5,rotate=90]	
\draw [fill=gray!45,gray!45] (-.6,-.6)  rectangle (.6,.6)   ;
\draw[line width=1.2] (-0.4,-0.52) ..controls +(.3,.5).. (-.4,.53);
\draw[line width=1.2] (0.4,-0.52) ..controls +(-.3,.5).. (.4,.53);
\end{tikzpicture}
\hspace{.5cm}
\text{ and }\hspace{.5cm}
\begin{tikzpicture}[baseline=-0.4ex,scale=0.5,rotate=90] 
\draw [fill=gray!45,gray!45] (-.6,-.6)  rectangle (.6,.6)   ;
\draw[line width=1.2,black] (0,0)  circle (.4)   ;
\end{tikzpicture}
= -(A^2+A^{-2}) 
\begin{tikzpicture}[baseline=-0.4ex,scale=0.5,rotate=90] 
\draw [fill=gray!45,gray!45] (-.6,-.6)  rectangle (.6,.6)   ;
\end{tikzpicture}
;
$$

$$
\begin{tikzpicture}[baseline=-0.4ex,scale=0.5,>=stealth]
\draw [fill=gray!45,gray!45] (-.7,-.75)  rectangle (.4,.75)   ;
\draw[->] (0.4,-0.75) to (.4,.75);
\draw[line width=1.2] (0.4,-0.3) to (0,-.3);
\draw[line width=1.2] (0.4,0.3) to (0,.3);
\draw[line width=1.1] (0,0) ++(90:.3) arc (90:270:.3);
\draw (0.65,0.3) node {\scriptsize{$+$}}; 
\draw (0.65,-0.3) node {\scriptsize{$+$}}; 
\end{tikzpicture}
=
\begin{tikzpicture}[baseline=-0.4ex,scale=0.5,>=stealth]
\draw [fill=gray!45,gray!45] (-.7,-.75)  rectangle (.4,.75)   ;
\draw[->] (0.4,-0.75) to (.4,.75);
\draw[line width=1.2] (0.4,-0.3) to (0,-.3);
\draw[line width=1.2] (0.4,0.3) to (0,.3);
\draw[line width=1.1] (0,0) ++(90:.3) arc (90:270:.3);
\draw (0.65,0.3) node {\scriptsize{$-$}}; 
\draw (0.65,-0.3) node {\scriptsize{$-$}}; 
\end{tikzpicture}
=0,
\hspace{.2cm}
\begin{tikzpicture}[baseline=-0.4ex,scale=0.5,>=stealth]
\draw [fill=gray!45,gray!45] (-.7,-.75)  rectangle (.4,.75)   ;
\draw[->] (0.4,-0.75) to (.4,.75);
\draw[line width=1.2] (0.4,-0.3) to (0,-.3);
\draw[line width=1.2] (0.4,0.3) to (0,.3);
\draw[line width=1.1] (0,0) ++(90:.3) arc (90:270:.3);
\draw (0.65,0.3) node {\scriptsize{$+$}}; 
\draw (0.65,-0.3) node {\scriptsize{$-$}}; 
\end{tikzpicture}
=A^{-1/2}
\begin{tikzpicture}[baseline=-0.4ex,scale=0.5,>=stealth]
\draw [fill=gray!45,gray!45] (-.7,-.75)  rectangle (.4,.75)   ;
\draw[-] (0.4,-0.75) to (.4,.75);
\end{tikzpicture}
\hspace{.1cm} \text{ and }
\hspace{.1cm}
A^{1/2}
\heightexch{->}{-}{+}
- A^{5/2}
\heightexch{->}{+}{-}
=
\heightcurve.
$$
The product of two classes of stated tangles $[T_1,s_1]$ and $[T_2,s_2]$ is defined by  isotoping $T_1$ and $T_2$  in $\Sigma \times (1/2, 1) $ and $\Sigma \times (0, 1/2)$ respectively and then setting $[T_1,s_1]\cdot [T_2,s_2]=[T_1\cup T_2, s_1\cup s_2]$.
\par We obtain two functors $\mathcal{S}_A: \MS \to \Mod_k$ and $\mathcal{S}_A: \MS^{str} \to \Alg_k$. 
\end{definition}

\underline{\textbf{Splitting morphisms:}} For $a,b$ two boundary edges of $\mathbf{\Sigma}$, there exists an injective morphism of algebras $\theta_{a\# b} : \mathcal{S}_A(\mathbf{\Sigma}_{a\# b}) \hookrightarrow \mathcal{S}_A(\mathbf{\Sigma})$, named the \textbf{splitting morphism}, defined as follows. Let $c \subset \Sigma_{a\# b}$ be the common image of $a$ and $b$ and think of $\Sigma$ as being obtained from $\Sigma_{a\# b}$ by splitting along $c$. For $(T,s)$ a stated tangle in $\mathbf{\Sigma}_{a\# b}$, isotope $T$ such that it intersects $\partial_c T:=c\times (0,1)$ transversally along points with pairwise distinct heights and framing in the $(0,1)$ direction towards $1$. Then split $T$ along $c\times (0,1)$ to obtain a tangle $T(c)$ in $\mathbf{\Sigma}$. Each intersection point $v\in \partial_c T$ gives rise to two points $v_a \in a$ and $v_b \in b$ of $\partial T(c)$ and we have a partition 
$\partial T(c)= \partial T  \bigsqcup_{v\in \partial_c T} \{v_a, v_b\}  $.
 A state $s': \partial T(c) \to \{-, +\}$ is \textit{admissible} if its restriction to $\partial T \subset \partial T(c)$ is equal to $s$ and for each $v \in \partial_c T$ then $s(v_a)=s(v_b)$. The splitting morphism is then defined by 
$$ \theta_{a\# b} : \mathcal{S}_A(\mathbf{\Sigma}_{a\# b}) \hookrightarrow \mathcal{S}_A(\mathbf{\Sigma}), \quad [T,s]\mapsto \sum_{s' \mathrm{ admissible}} [T(c), s'].$$
Definition \ref{def_sskein} is designed to ensure that $\theta_{a\# b} ([T,s])$ does not depend on the way we have isotoped $T$ before splitting and that $\theta_{a\# b}$ is a well-defined injective morphism of algebras (see \cite{LeStatedSkein}). 

\par \underline{\textbf{Bigon:}} The \textbf{bigon} $\mathbb{B}=(D^2, \{a_L, a_R\})=\adjustbox{valign=c}{\includegraphics[width=0.5cm]{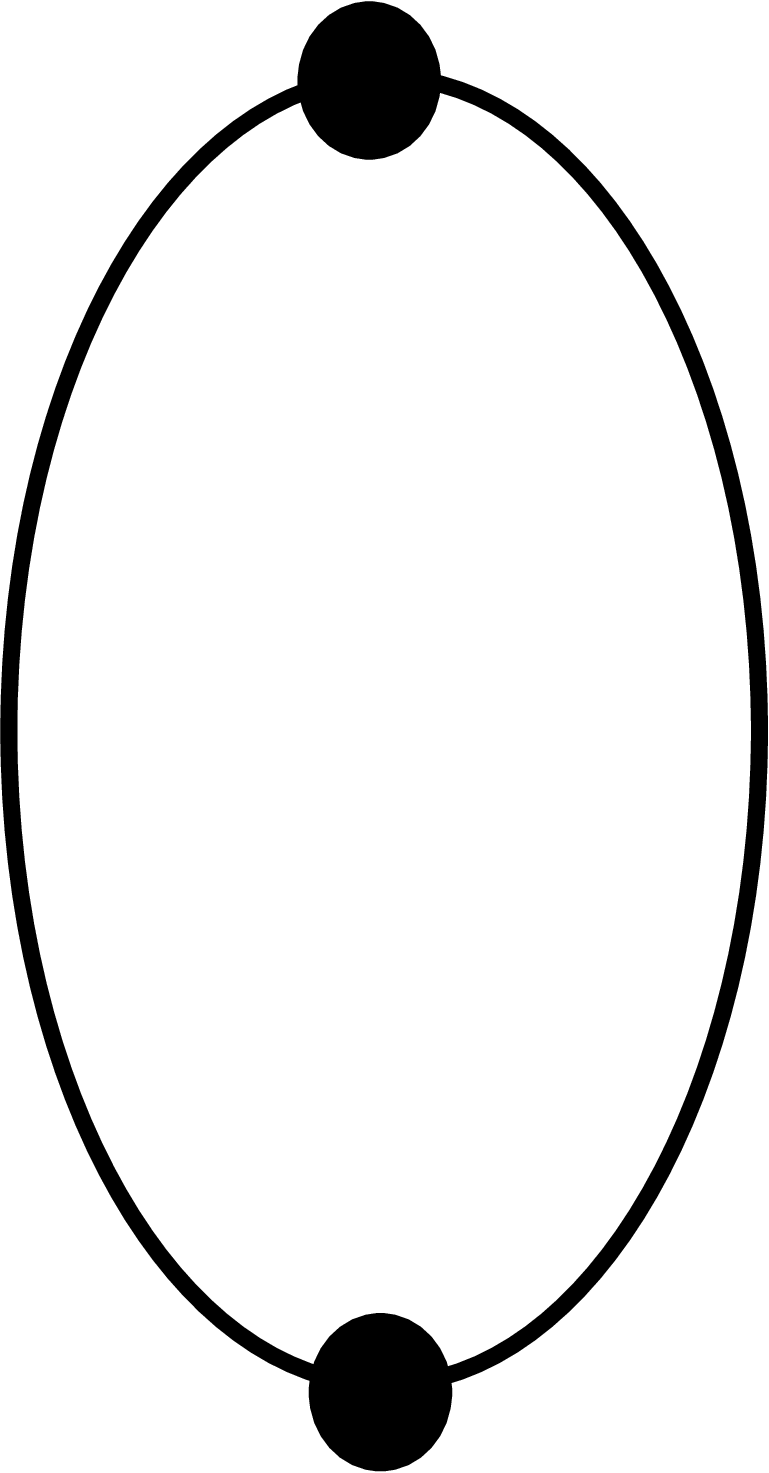}}$ is the marked surface made of a disc with two boundary edges $a_L, a_R$. When gluing two bigon $\mathbb{B}$ and $\mathbb{B}'$ together by gluing $a_R$ to $a_L'$ we obtain another bigon. The splitting morphism $\Delta:= \theta_{a_R \# a_L'} : \mathcal{S}_A(\mathbb{B})^{\otimes 2} \to \mathcal{S}_A(\mathbb{B})$ endows $\mathcal{S}_A(\mathbb{B})$ with a structure of Hopf algebra which were recognized to be isomorphic to the quantum group $\mathcal{O}_q \SL_2$ in \cite{KojuQuesneyClassicalShadows, CostantinoLe19} where $q:= A^2$.  

\par \underline{\textbf{Comodule structures:}} Let $a$ be a boundary edge of $\mathbf{\Sigma}$. While gluing $\mathbb{B}$ with $\mathbf{\Sigma}$ by gluing $a_R$ to $a$, we obtain a marked surface isomorphic to $\mathbf{\Sigma}$. Therefore the splitting morphism $\Delta_a^L := \theta_{a_R \# a} : \mathcal{S}_A(\mathbf{\Sigma}) \to \mathcal{O}_q\SL_2 \otimes \mathcal{S}_A(\mathbf{\Sigma})$ endows $\mathcal{S}_A(\mathbf{\Sigma})$ with a structure of left $\mathcal{O}_q\SL_2$ comodule. Similarly, while gluing $a$ with $a_L$ the splitting morphism $\Delta_a^R : \mathcal{S}_A(\mathbf{\Sigma}) \to \mathcal{S}_A(\mathbf{\Sigma}) \otimes \mathcal{O}_q\SL_2$ defines a right comodule structure. 

\begin{theorem}(Fundamental theorem of stated skein algebras \cite{KojuQuesneyClassicalShadows, CostantinoLe19})\label{theorem_splitting}:
\par For $a,b$ two boundary edges of $\mathbf{\Sigma}$, we have a left exact sequence: 
$$ 0 \to \mathcal{S}_A(\mathbf{\Sigma}_{a\# b}) \xrightarrow{\theta_{a\# b}} \mathcal{S}_A(\mathbf{\Sigma}) \xrightarrow{\Delta^L_a - \sigma \circ \Delta^R_b} \mathcal{O}_q \SL_2 \otimes \mathcal{S}_A(\mathbf{\Sigma}), $$
where $\sigma (x\otimes y):= y\otimes x$. 
\end{theorem}
In particular, $\mathcal{S}_A(\mathbf{\Sigma}_{a\# b})$ is completely determined by $\mathcal{S}_A(\mathbf{\Sigma})$ with its bicomodule structure.
\par \underline{\textbf{Triangulations:}}
The \textbf{triangle} $\mathbb{T}=(D^2, \{e_1, e_2, e_3\})=\adjustbox{valign=c}{\includegraphics[width=0.9cm]{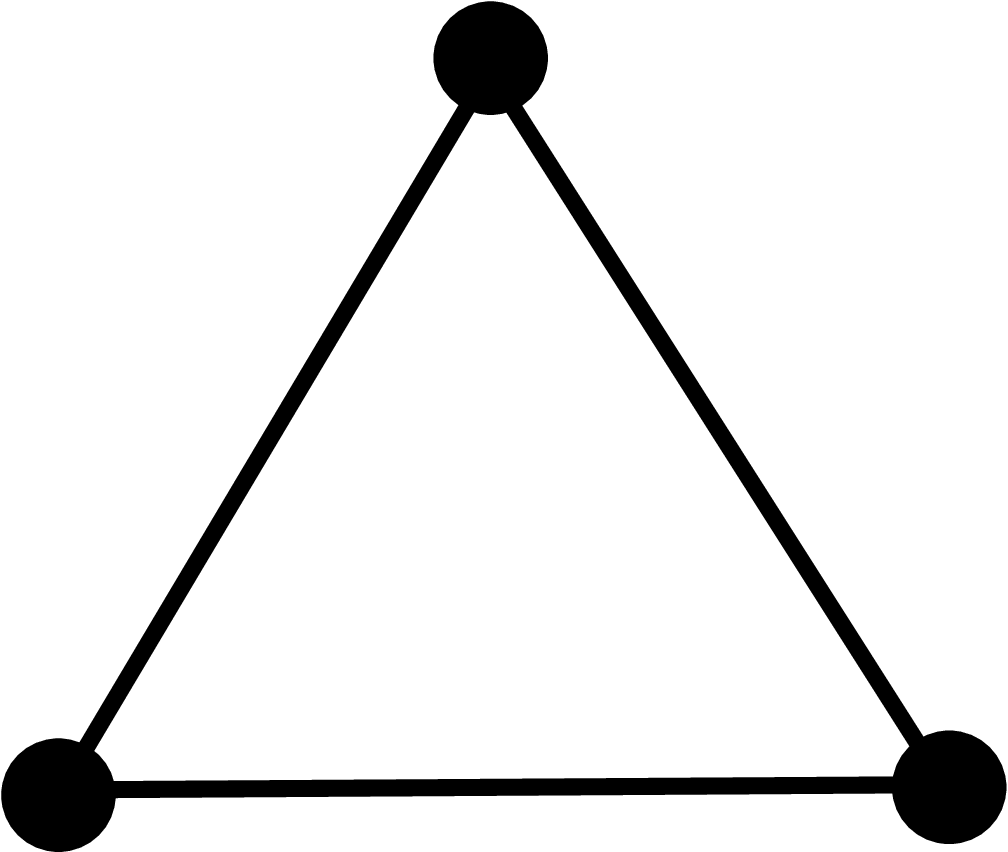}}$ is a disc with three boundary edges; we suppose that $e_1, e_2, e_3$ are cyclically ordered in the clockwise order. A marked surface $\mathbf{\Sigma}$ is \textbf{triangulable} if it can be obtained from a disjoint union of triangles by gluing some pairs of boundary edges; we then call \textbf{triangulation} $\Delta$  the combinatoric data made of the disjoint union of triangles (the faces with set $F(\Delta)$) and the pairs of glued edges. The splitting morphism defines an embedding $\theta^{\Delta} : \mathcal{S}_A(\mathbf{\Sigma}) \hookrightarrow \otimes_{\mathbb{T}\in F(\Delta)} \mathcal{S}_A(\mathbb{T})$ and Theorem \ref{theorem_splitting} implies that $\mathcal{S}_A(\mathbf{\Sigma}) $ is completely determined by the triangle algebra $\mathcal{S}_A(\mathbb{T})$ together with the combinatoric of the triangulation. This was the original motivation for the authors of \cite{BonahonWongqTrace} to introduce stated skein algebras. 

\par \underline{\textbf{Fusion and transmutation:}} Let $a,b$ be two boundary edges of $\mathbf{\Sigma}$. The \textbf{fusion of} $\mathbf{\Sigma}$ along $a$ and $b$ is the marked surface $\mathbf{\Sigma}_{a\circledast b}$ obtained by gluing $\mathbf{\Sigma}$ along a triangle $\mathbb{T}$ by gluing $a$ with $e_1$ and $b$ with $e_2$. Let $c$ be the boundary edge of $\mathbf{\Sigma}_{a\circledast b}$ corresponding to $e_3$. There is a natural (not strong) embedding $i: \mathbf{\Sigma} \hookrightarrow \mathbf{\Sigma}_{a\circledast b}$ which is the identity map outside some collar neighborhoods  of $e_1$ and $e_2$  and sends these collar neighborhoods into $\mathbb{T} \subset \mathbf{\Sigma}_{a\circledast b}$ by embedding $a$ and $b$ into $c$ with $a$ with higher height  than $b$. Since both $a$ and $b$ are sent into $c$, $i$ is not a strong embedding so the induced linear map $i_*: \mathcal{S}_A(\mathbf{\Sigma}) \to \mathcal{S}_A(\mathbf{\Sigma}_{a\circledast b})$ is not an algebra morphism.

\begin{theorem}(Fusion theorem \cite{CostantinoLe19}) \label{theorem_fusion}
\begin{enumerate}
\item 
$i_*: \mathcal{S}_A(\mathbf{\Sigma}) \to \mathcal{S}_A(\mathbf{\Sigma}_{a\circledast b})$ is a linear isomorphism.
\item  If $\mu$ denotes the product of $\mathcal{S}_A(\mathbf{\Sigma})$ and  $\mu_{a\circledast b}$ denotes the pull-back of $\mu$ by $i_*$, then $\mu_{a\circledast b}$ factorizes as 
\begin{equation}\label{eq_fusion} \mu_{a \circledast b}: \mathcal{S}_A(\mathbf{\Sigma})^{\otimes 2}  \xrightarrow{\Delta^L_a \otimes \Delta^L_b} \mathcal{O}_q\SL_2 \otimes \mathcal{S}_A(\mathbf{\Sigma})\otimes  \mathcal{O}_q\SL_2 \otimes \mathcal{S}_A(\mathbf{\Sigma}) \xrightarrow{\id \otimes \sigma \otimes \id}  (\mathcal{O}_q\SL_2)^{\otimes 2} \otimes \mathcal{S}_A(\mathbf{\Sigma})^{\otimes 2} \xrightarrow{r \otimes \mu} \mathcal{S}_A(\mathbf{\Sigma}),\end{equation}
where $\sigma$ is the flip and $r$ is the co-R matrix of the coribbon Hopf algebra $\mathcal{O}_q\SL_2$. 
\item The pull-back of the comodule map $\Delta_c^L$ by $i_*$ is the composition
$$ \mathcal{S}_A(\mathbf{\Sigma}) \xrightarrow{ ( \Delta_b^L \otimes \id )\Delta_a^L} (\mathcal{O}_q\SL_2) ^{\otimes 2} \otimes \mathcal{S}_A(\mathbf{\Sigma}) \xrightarrow{ \mu \otimes \id}  \mathcal{O}_q\SL_2 \otimes \mathcal{S}_A(\mathbf{\Sigma}).$$
\end{enumerate}
\end{theorem}
Let us stress two applications of this theorem. When fusioning the two boundary edges of the bigon, we obtain a once-punctured monogon $\mathbb{B}_{a_L \circledast a_R}= \mathbf{m}_1$ which consists in an annulus with a single boundary edge. So as a vector space $\mathcal{S}_A(\mathbf{m}_1) \cong \mathcal{O}_q\SL_2$ with a product $\mu_{a_L\circledast a_R}$ and a comodule structure given by the formula in Theorem \ref{theorem_fusion}. The coproduct of $\mathcal{O}_q\SL_2$ endows $\mathcal{S}_A(\mathbf{m_1})$ with the structure of a (braided) Hopf algebra object in the braided category $\mathcal{O}_q\SL_2-\Comod$ named the \textbf{braided quantum group} by Majid (\cite{Majid_QGroups}) and denoted $B_q\SL_2$. The processes from passing from the co-braided Hopf algebra $\mathcal{O}_q\SL_2$ to its braided version $B_q\SL_2$ by changing the product using Equation \eqref{eq_fusion} and comodule structure is called the \textbf{transmutation}. 
\par In the particular case where $\mathbf{\Sigma}=\mathbf{\Sigma}_1 \sqcup \mathbf{\Sigma}_2$ with $a$ and $b$ boundary edges of $\mathbf{\Sigma}_1$ and $\mathbf{\Sigma}_2$, we write $\mathbf{\Sigma}_1 \circledast \mathbf{\Sigma}_2:= (\mathbf{\Sigma})_{a\circledast b}$. In this case, the algebra structure given by Equation \eqref{eq_fusion} is called the \textbf{braided tensor product} and we use the notation $\mathcal{S}_A(\mathbf{\Sigma}_1)\overline{\otimes} \mathcal{S}_A(\mathbf{\Sigma}_2)$ to denote this tensor product.
 Denote by $\mathbf{\Sigma}_{g,n}^*=(\Sigma_{g, n+1}, \{a\})$ the genus $g$ surface with $n+1$ boundary. Then we have $\mathbf{\Sigma}_{g_1, n_1}^* \circledast \mathbf{\Sigma}_{g_2, n_2}^* = \mathbf{\Sigma}_{g_1+g_2, n_1+n_2}^* $ and $\mathbf{\Sigma}_{0,1}^* =\mathbf{m}_1$ so the theorem implies 
 $$ \mathcal{S}_A(\mathbf{\Sigma}_{g,n}^*) \cong \mathcal{S}_A(\mathbf{\Sigma}_{1,0}^*)^{\overline{\otimes} g} \overline{\otimes} (B_q\SL_2)^{\overline{\otimes} n}.$$
 \par \underline{\textbf{Presenting graphs:}} Call a marked surface $\mathbf{\Sigma}=(\Sigma, \mathcal{A})$ \textbf{essential} if $\Sigma$ is connected and $\mathcal{A}$ not empty. In this case, for each boundary edge $a\in \mathcal{A}$ fix a base point $v_a \in a$ and let $\mathbb{V}= \{ v_a\}_{a\in \mathcal{A}}$ be the set of such base points. A \textbf{presenting graph} for $\mathbf{\Sigma}$ is an embedded oriented graph $\Gamma \subset \Sigma$ whose set of vertices is $\mathbb{V}$ and such that $\Sigma$ retracts on $\Gamma$. Denote by $E(\Gamma)$ the set of edges. A simple argument detailed in \cite{CostantinoLe19} shows that $\mathbf{\Sigma}$ can be obtained from a disjoint union of $|E(\Gamma)|$ bigons by fusioning some pairs of boundary edges together so Theorem \ref{theorem_fusion} implies that we have an isomorphism of vector spaces 
 $$\Psi^{\Gamma}:  \mathcal{S}_A(\mathbf{\Sigma}) \cong (\mathcal{O}_q\SL_2)^{\otimes E(\Gamma)}.$$
\par \textbf{\underline{Basic properties of stated skein algebras:}} 
\begin{itemize}
\item $\mathcal{S}_A(\mathbf{\Sigma})$ \textit{is free over the ground ring} $k$: for unmarked surfaces, a simple classical argument shows that the set of classes of multicurves without contractible component forms a basis of $\mathcal{S}_A(\Sigma)$. For essential marked surface, we can use the above isomorphism $\Psi^{\Gamma}$ and the fact that $\mathcal{O}_q\SL_2$ is free (see \cite{LeStatedSkein, KojuKaruo_RepRSSkein} for alternative bases). 
\item $\mathcal{S}_A(\mathbf{\Sigma})$ \textit{is a domain}: for the bigon, it is well known that $\mathcal{O}_q\SL_2$ is a domain. Using Theorem \ref{theorem_fusion}, we deduce that $\mathcal{S}_A(\mathbb{T})\cong \mathcal{O}_q\SL_2 \overline{\otimes} \mathcal{O}_q\SL_2$ is a domain too and using the embedding $\theta^{\Delta} : \mathcal{S}_A(\mathbf{\Sigma}) \hookrightarrow \otimes_{\mathbb{T} \in F(\Delta)} \mathcal{S}_A(\mathbf{\Sigma})$ we obtain the result for any triangulable marked surfaces. The proof for closed surfaces is based on graduation and is more technical (see \cite{PrzytyckiSikora_SkeinDomain}).
\item $\mathcal{S}_A(\mathbf{\Sigma})$ \textit{if finitely generated}; \textit{when} $\mathbf{\Sigma}$ \textit{is essential, it is even finitely presented}: for unmarked surfaces, the result has been proved in several ways (\cite{BullockGeneratorsSkein, AbdielFrohman_SkeinFrobenius, FrohmanKania_SkeinRootUnity, SantharoubaneSkeinGenerators}).
It is well known that  $\mathcal{O}_q\SL_2$ if finitely generated by the four matrix elements functions $a,b,c,d$ of the standard representation of $U_q\mathfrak{sl}_2$. When $\mathbf{\Sigma}$ is essential, a set of generators for $\mathcal{S}_A(\mathbf{\Sigma})$ is given by the set of elements $x_e:= (\Psi^{\Gamma})^{-1}(1\otimes \ldots \otimes 1 \otimes  x\otimes 1 \otimes \ldots \otimes 1)$ ($x$ is in the $e$-th factor) with $x \in \{a,b,c,d\}$ and $e\in E(\Gamma)$ (\cite{KojuAzumayaSkein}). Using skein (height exchange) relations, for $e,e'\in E(\Gamma)$ and $x, x' \in \{a,b,c,d\}$ we can write the product $x_e x'_{e'}$ as a linear combination of monomials of the form $(x_i)_{e'}(x'_i)_e$ thus obtaining a relation $x_ex'_{e'} = \sum_i \alpha_i (x_i)_{e'}(x'_i)_e$. Together with the relations $.$ for $e\in E(\Gamma)$ these relations define a finite presentation of $\mathcal{S}_A(\mathbf{\Sigma})$ (\cite{KojuPresentationSSkein}). In the particular case where $\mathbf{\Sigma}=\mathbf{\Sigma}_{g,n}^*$ and $\Gamma$ is the daisy graph, the obtained presentation of $\mathcal{S}_A(\mathbf{\Sigma}_{g,n}^*)$ is precisely the presentation defining the lattice gauge field theoretical algebra $\mathcal{L}_{g,n}$ from Matthieu's talk. We thus obtain the following theorem already explained by Matthieu.

\end{itemize}

\begin{theorem}(\cite{Faitg_LGFT_SSkein}) We have an isomorphism of algebras $\mathcal{S}_A(\mathbf{\Sigma}_{g,n}^*) \cong \mathcal{L}_{g,n}$. \end{theorem}

\section{Semi classical limits of stated skein algebras}

When $A^{1/2}=+1$,  $\mathcal{S}_{+1}(\mathbf{\Sigma})$ is commutative and the scheme $X(\mathbf{\Sigma}):= \Specm\left(\mathcal{S}_{+1}(\mathbf{\Sigma})\right)$ is a variety (over $\mathbb{C}$) with a Poisson structure described as follows. Consider the ring $k_{\hbar}:=\mathbb{C}[[\hbar]]$ of formal power series and $A^{1/2}:= \exp(\hbar/2)\in k_{\hbar}$ and denote by $\mathcal{S}_{\hbar}(\mathbf{\Sigma})$ the corresponding stated skein algebra. Since stated skein algebras are free, we have a natural linear isomorphism $\mathcal{S}_{+1}(\mathbf{\Sigma})\otimes_{\mathbb{C}}k_{\hbar} \cong \mathcal{S}_{\hbar}(\mathbf{\Sigma})$; denote by $\star$ the pull-back of the product of $\mathcal{S}_{\hbar}(\mathbf{\Sigma})$ by this isomorphism. The Poisson structure on $X(\mathbf{\Sigma})$  is then defined by the formula: 
$$ x \star y - y\star x \equiv \hbar \{x,y\} \pmod{\hbar^2} \quad \mbox{, for }x, y \in \mathcal{S}_{+1}(\mathbf{\Sigma}).$$
$X(\mathbf{\Sigma})$ admits  the following geometric interpretation. 
Suppose  $\mathbf{\Sigma}$ essential. The \textbf{fundamental groupoid} $\pi_1(\Sigma, \mathbb{V})$ with base points $\mathbb{V}$ is the groupoid whose set of objects is $\mathbb{V}$ and whose morphisms $\alpha: v_1 \to v_2$ are homotopy classes of paths $c_{\alpha}: [0,1] \to \Sigma$ with $s(\alpha):=c(0)=v_1$ and $t(\alpha):=c(1)=v_2$. Composition is the concatenation of paths. A \textbf{representation} is a functor $\rho: \pi_1(\Sigma, \mathbb{V})\to \SL_2$ where $\SL_2$ is seen as a category with only one element whose endomorphism group is $\SL_2$. When $\mathbb{V}=\{v\}$ has a single element, this is a same as a group morphism $\rho: \pi_1(\Sigma, v) \to \SL_2$ in the usual sense.
 The set $\mathcal{R}_{\SL_2}(\mathbf{\Sigma}):=\Hom(\pi_1(\Sigma, \mathbb{V}), \SL_2)$ has a natural structure of affine variety and a Poisson structure defined in \cite{FockRosly}. When $\mathbb{V}=\{v\}$ has one element, $\SL_2$ acts algebraically of $\Hom(\pi_1(\Sigma, v), \SL_2)$ and the GIT quotient $\mathcal{X}(\Sigma):= \Hom(\pi_1(\Sigma, v), \SL_2)\sslash \SL_2$ is called the \textbf{character variety}. 
 
 \begin{theorem}(Classical limits of stated skein algebras) \label{theorem_classical_limit}
 \begin{enumerate}
 \item \cite{Bullock, PS00, Turaev91} If $\mathbf{\Sigma}$ is unmarked, we have a Poisson isomorphism $X(\mathbf{\Sigma}) \cong \mathcal{X}(\Sigma)$ where the character variety $\mathcal{X}(\Sigma)$ is equipped with its Atiyah-Bott-Goldman Poisson structure from \cite{AB, Goldman86}. 
 \item \cite{KojuQuesneyClassicalShadows} If $\mathbf{\Sigma}$ is essential, we have a Poisson isomorphism $X(\mathbf{\Sigma})\cong \mathcal{R}_{\SL_2}(\Sigma, \mathbb{V})$ where the representation variety $\mathcal{R}_{\SL_2}(\mathbf{\Sigma})$ is equipped with its Fock-Rosly Poisson structure from \cite{FockRosly}.
 \end{enumerate}
 \end{theorem}
 
 Let us introduce a partition of character varieties. Consider the injective group morphisms $\mathbb{Z}/2\mathbb{Z} \hookrightarrow \mathbb{C}^* \hookrightarrow{\SL_2}$; the first has image $\{-1, +1\}\subset \mathbb{C}^*$ and the second sends $z$ to the matrix $\begin{pmatrix} z & 0 \\ 0 & z^{-1} \end{pmatrix}$. They induce some embeddings $\mathrm{H}^1(\Sigma; \mathbb{Z}/2\mathbb{Z}) \hookrightarrow \mathrm{H}^1(\Sigma; \mathbb{C}^*) \hookrightarrow \mathcal{X}(\Sigma)$. A class in the image $\mathrm{H}^1(\Sigma; \mathbb{Z}/2\mathbb{Z}) \hookrightarrow \mathcal{X}(\Sigma)$ is called a \textbf{central representation} whereas a class in the image $\mathrm{H}^1(\Sigma; \mathbb{C}^*) \hookrightarrow \mathcal{X}(\Sigma)$ which is not central is called a \textbf{diagonal representation}. A class which is neither central, nor diagonal is the class of an irreducible representation we thus obtain a partition
 $$ \mathcal{X}(\Sigma)= \mathcal{X}^{central}(\Sigma) \sqcup \mathcal{X}^{diagonal}(\Sigma) \sqcup \mathcal{X}^{irred}(\Sigma).$$
 
 \section{Representations of stated skein algebras at roots of unity}
 
 \underline{\textbf{Frobenius morphisms:}} Let $\zeta^{1/2}\in \mathbb{C}$ be a root of unity of odd order $N\geq 3$ and consider the stated skein algebra $\mathcal{S}_{\zeta}(\mathbf{\Sigma})$ where $A=\zeta$. 
 
 \begin{theorem}(\cite{BonahonWong1} for unmarked surfaces, \cite{KojuQuesneyClassicalShadows} in general) There exists an injective  morphism of algebras $Fr_{\mathbf{\Sigma}}: \mathcal{S}_{+1}(\mathbf{\Sigma}) \hookrightarrow \mathcal{S}_{\zeta}(\mathbf{\Sigma})$, named \textbf{Frobenius morphism}, whose image lies in the center of $\mathcal{S}_{\zeta}(\mathbf{\Sigma})$. 
 \end{theorem}
 
 An irreducible representation $r: \mathcal{S}_{\zeta}(\mathbf{\Sigma}) \to \End(V)$ sends central elements to scalars, so it induces a character over $\mathcal{S}_{+1}(\mathbf{\Sigma})\cong \mathcal{O}[X(\mathbf{\Sigma})]$ and therefore a point in $X(\mathbf{\Sigma})$ named its \textbf{classical shadow}.
 
 \par \underline{\textbf{Reduced stated skein algebras:}} Let $p$ be a boundary puncture of $\mathbf{\Sigma}$. Let $\alpha(p)$ be the corner arc around $p$ and $\alpha(p)_{-+}$ be the arc $\alpha(p)$ with states $-$ and $+$, i.e. $\alpha(p)_{-+}=\adjustbox{valign=c}{\includegraphics[width=1.9cm]{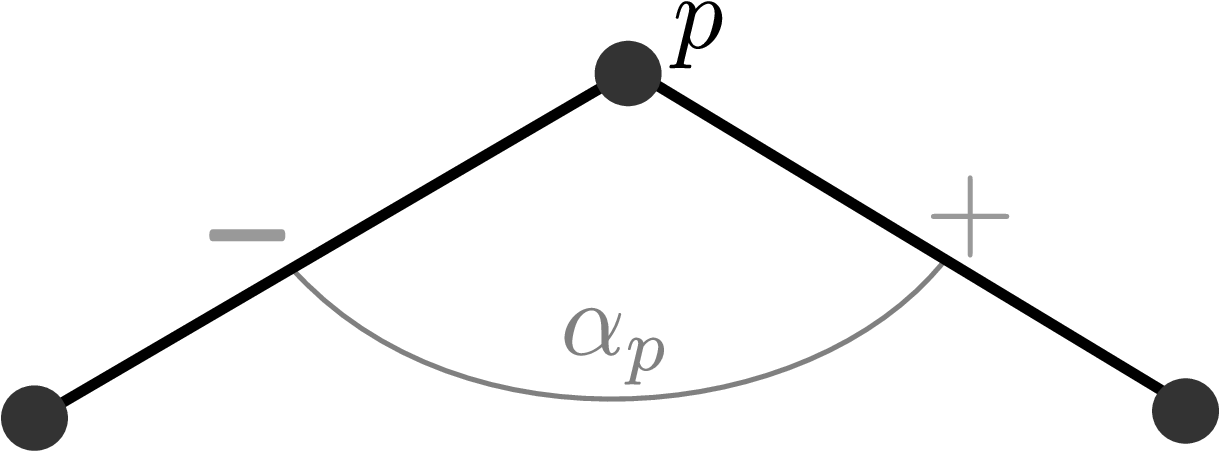}}$.  We call $\alpha(p)_{-+} \in \mathcal{S}_A(\mathbf{\Sigma})$ a \textbf{bad arc}. The \textbf{reduced stated skein algebra} $\overline{\mathcal{S}}_A(\mathbf{\Sigma})$ is the quotient of $\mathcal{S}_A(\mathbf{\Sigma})$ by the ideal generated by bad arcs. Most of properties of stated skein algebras still hold for their reduced version: $\overline{\mathcal{S}}_A(\mathbf{\Sigma})$ is also free, is a domain, has finite presentation, the splitting morphisms and Frobenius morphisms are well-defined in the reduced case by passing to the quotient. However there is no exact sequence of the form of the one in Theorem \ref{theorem_splitting}, i.e. $\overline{\mathcal{S}}_A(\mathbf{\Sigma}_{a\# b})$ is not determined by $\overline{\mathcal{S}}_A(\mathbf{\Sigma})$ and Theorem \ref{theorem_fusion} has no analogue. For an essential marked surface $\mathbf{\Sigma}$ then $\overline{X}(\mathbf{\Sigma}):=\Specm(\overline{\mathcal{S}}_{+1}(\mathbf{\Sigma}))$ is isomorphic to the subvariety of $\mathcal{R}_{\SL_2}(\mathbf{\Sigma})$ of representations $\rho: \pi_1(\Sigma, \mathbb{V})\to \SL_2$ which send each corner arc $\alpha(p)$ to an element of the small Bruhat cell of $\SL_2$, i.e. to $\SL_2$ matrices with vanishing  top left matrix coefficient (see \cite{KojuKaruo_RepRSSkein}).

 \par  \underline{\textbf{Center:}} In general, the centers $Z_{\mathbf{\Sigma}}$ and $\overline{Z}_{\mathbf{\Sigma}}$ of $\mathcal{S}_{\zeta}(\mathbf{\Sigma})$ and $\overline{\mathcal{S}}_{\zeta}(\mathbf{\Sigma})$ are larger than the "small centers"  $Z^0_{\mathbf{\Sigma}}\subset Z_{\mathbf{\Sigma}}$ and $\overline{Z}^0_{\mathbf{\Sigma}}\subset \overline{Z}_{\mathbf{\Sigma}}$ defined as being the images of the Frobenius morphisms. For $p$ an inner puncture of $\mathbf{\Sigma}$ then the peripheral curve $\gamma_p$ which encircles $p$ once defines a central element in both $\mathcal{S}_{\zeta}(\mathbf{\Sigma})$ and $\overline{\mathcal{S}}_{\zeta}(\mathbf{\Sigma})$. Also for $\partial$ a boundary component $\mathbf{\Sigma}$ which contains some boundary edges, then the element $\alpha_{\partial}\in \overline{\mathcal{S}}_{\zeta}(\mathbf{\Sigma})$,  defined as the product of stated corner arcs $\alpha(p)_{++}$ for $p\subset \partial$, is also central and admits an inverse $\alpha_{\partial}^{-1}$ which is the product of elements $\alpha(p)_{--}$ for  $p\subset \partial$. If $\partial$ is a boundary component of $\Sigma$ which contains an even positive number of boundary edges, we can define a central element $\beta_{\partial} \in \mathcal{S}_{\zeta}(\mathbf{\Sigma})$ (see \cite{Yu_CenterSSKein}). 
 
 \begin{theorem}(Centers of stated skein algebras \cite{FrohmanKaniaLe_UnicityRep, KojuAzumayaSkein, GanevJordanSafranov_FrobeniusMorphism, KojuKaruo_RepRSSkein, Yu_CenterSSKein})
 \begin{enumerate}
 \item The center $\overline{Z}_{\mathbf{\Sigma}}$ of $\overline{\mathcal{S}}_{\zeta}(\mathbf{\Sigma})$ is generated by $\overline{Z}^0_{\mathbf{\Sigma}}$, the peripheral curves $\gamma_p$ and the boundary central elements $\alpha_{\partial}^{\pm 1}$. Moreover  we have a bijection
 $$ \widehat{\overline{X}}(\mathbf{\Sigma}):= \Specm(\overline{Z}_{\mathbf{\Sigma}})\cong \{ (\rho, h_{p_1}, \ldots, h_{p_n}, h_{\partial_1}, \ldots, h_{\partial_m}):  \rho \in \overline{X}(\mathbf{\Sigma}), T_N(h_{p_i})=-\tr(\rho(\gamma_{p_i})), h_{\partial_j}^N = \rho(\alpha_{\partial})_{++} \}. $$
 Here $T_N(X)$ is the Chebyshev polynomial of first kind and $\rho$ is seen either as a representation $\rho: \pi_1(\mathbf{\Sigma}, \mathbb{V})\to \SL_2$ when $\mathbf{\Sigma}$ is essential or as a class in the GIT quotient $\mathcal{X}(\mathbf{\Sigma})$ in the unmarked case. 
 \item The center $Z_{\mathbf{\Sigma}}$ of $\mathcal{S}_{\zeta}(\mathbf{\Sigma})$ is generated by $Z^0_{\mathbf{\Sigma}}$ together with the peripheral curves $\gamma_p$ associated to inner punctures and the central elements $\beta_{\partial}$ associated to even marked boundary components. For $\mathbf{\Sigma}_{g,n}^*$, we have a bijection: 
  $$ \widehat{X}(\mathbf{\Sigma}_{g,n}^*):= \Specm({Z}_{\mathbf{\Sigma}})\cong \{ (\rho, h_{p_1}, \ldots, h_{p_n}):  \rho \in \Hom(\pi_1(\Sigma_{g,n+1}, v), \SL_2), T_N(h_{p_i})=-\tr(\rho(\gamma_{p_i}))\}. $$
 \end{enumerate}
 \end{theorem}

\par   \underline{\textbf{Azumaya loci:}}
 
 An irreducible representation $r: \mathcal{S}_{\zeta}(\mathbf{\Sigma}) \to \End(V)$ sends central elements to scalar operators therefore induces a character over $Z_{\mathbf{\Sigma}}$ so a point in $\widehat{X}(\mathbf{\Sigma})$ named its \textbf{full shadow}. An important property, first conjectured by Bonahon and Wong, is the fact that the character maps 
 $$ \chi: \Irrep\left( \mathcal{S}_{\zeta}(\mathbf{\Sigma}) \right) \to \widehat{X}(\mathbf{\Sigma}), \quad \chi: \Irrep\left( \overline{\mathcal{S}}_{\zeta}(\mathbf{\Sigma}) \right) \to \widehat{\overline{X}}(\mathbf{\Sigma}),$$
 are "almost" bijections. More precisely, there exists an open dense subset $\mathcal{AL} \subset \widehat{X}(\mathbf{\Sigma})$ (the Azumaya locus) such that the restriction $\chi : \chi^{-1}(\mathcal{AL}) \to \mathcal{AL}$ is a bijection: so a "generic" irreducible representation of $\mathcal{S}_{\zeta}(\mathbf{\Sigma})$ is completely determined by its full shadow (and similarly in the reduced case). Moreover, every such generic irreducible representations have the same dimension. 
 \par This property results from deep results of De Concini-Kac which we briefly summarized here. An \textbf{almost Azumaya algebra} is a complex algebra  $\mathscr{A}$ such that $(i)$ $\mathscr{A}$ is finitely generated as an algebra, $(ii)$ $\mathscr{A}$ is prime and $(iii)$ $\mathscr{A}$ is finitely generated as a module over its center $Z$.  Let $Q(Z)$ denote the fraction field of $Z$ obtained by localizing by every non zero element. By a theorem of Posner-Formaneck, $\mathscr{A}\otimes_Z Q(Z)$ is a central simple algebra and there exists a finite extension $F$ of $Q(Z)$ such that $\mathscr{A} \otimes_Z F\cong \Mat_D(F)$ is a matrix algebra of size $D$. The integer $D$ is called the \textit{PI-degree} of $\mathscr{A}$ and is characterized by the formula $\dim_{Q(Z)}(\mathscr{A}\otimes_ZQ(Z)) = D^2$. Let ${X}:= \Specm(Z)$ and for $x\in X$ corresponding to a maximal ideal $\mathfrak{m}_x$, consider the finite dimensional algebra $\mathscr{A}_x:= \quotient{\mathscr{A}}{\mathfrak{m}_x \mathscr{A}}$. 
  \par  The \textbf{Azumaya locus} of $\mathscr{A}$ is the subset
  $$ \mathcal{AL}:= \{x \in X : \mathscr{A}_x \cong \Mat_D(\mathbb{C}) \}.$$
  Since $\Mat_D(\mathbb{C})$ is simple with a single irreducible representation, for $x \in \mathcal{AL}$ there exists a unique irreducible representation (up to isomorphism) which full shadow $x$. We call this representation the \textbf{Azumaya representation with full shadow } $x$. 
  A deep theorem of De Concini Kac (\cite{DeConciniKacRepQGroups}, \cite[Theorem III.I.6-7]{BrownGoodearl}) asserts that the Azumaya locus of an almost Azumaya algebra is open dense. Moreover, any irreducible representation with full shadow outside the Azumaya locus has dimension $<D$ so the Azumaya locus admits the following characterization:

 $$\mathcal{AL}=\{ x \in \mathcal{X} | x \mbox{ is the shadow of an irreducible representation of maximal dimension }D\}.$$
 We already saw that $\mathcal{S}_{\zeta}(\mathbf{\Sigma})$ and $\overline{\mathcal{S}}_{\zeta}(\mathbf{\Sigma})$ are finitely generated and domain (thus prime). That they are finitely generated over their center easily results from the existence of the Frobenius morphism. 
 
 \begin{theorem}(Azumaya loci of skein algebras)\label{theorem_AL} Let $\mathbf{\Sigma}=(\Sigma_{g,n}, \mathcal{A})$ be a connected marked surface of genus $g$ with $n$ boundary component with $\chi(\Sigma_{g,n})<0$.
 \begin{enumerate}
 \item (\cite{FrohmanKaniaLe_UnicityRep} unmarked case, \cite{KojuAzumayaSkein} essential case) $\mathcal{S}_{\zeta}(\mathbf{\Sigma})$ and $\overline{\mathcal{S}}_{\zeta}(\mathbf{\Sigma})$  are almost Azumaya so their Azumaya loci $\mathcal{AL}(\mathbf{\Sigma})$ and $\overline{\mathcal{AL}}(\mathbf{\Sigma})$ are open dense.
 \item (\cite{FrohmanKaniaLe_UnicityRep, FrohmanKaniaLe_DimSkein} unmarked case, \cite{KojuAzumayaSkein} essential case) The PI-degree of $\overline{\mathcal{S}}_{\zeta}(\mathbf{\Sigma})$  is $N^{3g-3+n+|\mathcal{A}|}$.
 \item (\cite{GanevJordanSafranov_FrobeniusMorphism} for $\mathbf{\Sigma}_{g,0} ^*$, \cite{KojuKaruo_RepRSSkein} for $\mathbf{\Sigma}_{g,n}^*$, \cite{Yu_CenterSSKein} in general) The PI degree of $\mathcal{S}_{\zeta}(\mathbf{\Sigma})$ is $N^{3g-3+n_{ev}+\frac{3}{2}(n_{odd}+|\mathcal{A}|)}$ where $n_{ev}, n_{odd}$ are the number of boundary components having an even/odd number of boundary edges.
 \item (\cite{KojuKaruo_Azumaya}) For $\mathbf{\Sigma}=(\Sigma_{g,0}, \emptyset)$ an unmarked closed surface, the Azumaya locus $\mathcal{AL}(\Sigma)$ of $\mathcal{S}_{\zeta}(\Sigma_{g,0})$ is the complementary of the locus $\mathcal{X}^{central}$ of central representations. 
 \item Suppose $\mathbf{\Sigma}=(\Sigma_{g,n}, \emptyset)$ is unmarked, $n\geq 1$ and let $\widehat{x}=([\rho], h_{p_1}, \ldots, h_{p_n}) \in \widehat{X}(\mathbf{\Sigma})$. 
 \begin{itemize}
 \item (\cite{FKL_GeometricSkein}) $(i)$ If $g=0$, $[\rho]$ is irreducible and for all $i$, $\tr(\rho(\gamma_{p_i}))\neq \pm 2$, then $\widehat{x}\in \mathcal{AL}(\mathbf{\Sigma})$. $(ii)$ If for all $i$, $\tr(\rho(\gamma_{p_i}))\neq \pm 2$ and if $\tr(\rho(\gamma_{p_i}))=t_i +t_i^{-1}$ such that $t_1t_2\ldots t_n \neq 1$ then $\widehat{x}\in \mathcal{AL}(\mathbf{\Sigma})$. 
 \item (\cite{KojuKaruo_Azumaya}) If $[\rho]$ is central and  for all $i$, $h_{p_i}= -q^{n_i}-q^{-n_i}$ for some $n_i \in \{0, \ldots, N-1\}$ then $\widehat{x} \notin \mathcal{AL}(\mathbf{\Sigma})$. If $g\geq 1$, $[\rho]$ is diagonal not central and for all $i$, either $\tr(\rho(\gamma_{p_i}))\neq \pm 2$ or $h_{p_i}=\pm 2$, then $\widehat{x} \in \mathcal{AL}(\mathbf{\Sigma})$. 
 \item (\cite{KojuKaruo_Azumaya, Yu_CenterSSKein}) When $\Sigma= \Sigma_{1,1}$ or $\Sigma_{0,4}$, the Azumaya loci $\mathcal{AL}(\mathbf{\Sigma})$ are well known. 
 \end{itemize}
 \item (\cite{KojuKaruo_RepRSSkein}) Suppose $\mathbf{\Sigma}$ is essential and either has a boundary component with at least two boundary edges or which does not have any inner puncture. Then $\widehat{x}=(\rho, h_{p_i}, h_{\partial_j}) \in \overline{\mathcal{AL}}(\mathbf{\Sigma})$ if and only if for all inner puncture $p_i$, either $\rho(\gamma_{p_i})\neq \pm \mathds{1}_2$ or $h_{p_i}=\pm 2$.  
  \item (\cite{KojuKaruo_RepRSSkein}) The Azumaya locus of  $\mathcal{S}_{\zeta}(\mathbf{\Sigma}_{g,n}^*)\cong \mathcal{L}_{g,n}$ is the locus of elements $\widehat{x}=(\rho, h_{p_1}, \ldots, h_{p_n})$ such that $(1)$ for $p_{\partial}$ the unique boundary puncture then $\rho(\alpha(p_{\partial}))$ belongs to the big Bruhat cell of $\SL_2$ (its left top matrix coefficient does not vanishes) and $(2)$ for all $1\leq i \leq n$, either $\rho(\gamma_{p_i})\neq \pm \mathds{1}_2$ or $h_{p_i}=\pm 2$. 
 \end{enumerate}
 \end{theorem}
 
 \section{Mapping class group representations}
 
 \subsection{Deriving mapping class group representations from stated skein representations}
 {.}
 \\ \underline{\textbf{From skein representation to MCG representations:}} The mapping class group $\Mod(\Sigma)$ acts on $\mathcal{S}_{\zeta}(\mathbf{\Sigma})$ on the right as follows. For $\phi: \Sigma \to \Sigma$ an oriented homeomorphism, we denote by $\widetilde{\phi}: \Sigma\times [0,1]\to \Sigma \times [0,1]$ the homeomorphism defined by $\widetilde{\phi}(x,t):= (\phi(x), t)$. The right action of  $\Mod(\Sigma)$  on $\mathcal{S}_{\zeta}(\mathbf{\Sigma})$ is then defined by $[T,s]\cdot \phi:= [\widetilde{\phi}^{-1}(T), s\circ\widetilde{\phi}]$. This action induces  left actions of $\Mod(\Sigma)$ on $\widehat{X}(\mathbf{\Sigma})$ and $X(\mathbf{\Sigma})$. 
 \par 
 For  $G\subset \Mod(\Sigma)$ a subgroup, an irreducible representation $r: \mathcal{S}_{\zeta}(\mathbf{\Sigma}) \to \End(V)$ will be called \textbf{fixed by} $G$ if for every $\phi\in G$, the representation $\phi\cdot r : \mathcal{S}_{\zeta}(\mathbf{\Sigma}) \to \End(V), x \mapsto r( \phi^{-1}(x))$ is isomorphic to $r$. In this case, there exists an operator $L_r(\phi) \in \End(V)$, unique up to multiplication by a non zero scalar, such that 
$$ r(\phi^{-1}(x)) = L_r(\phi) r(x) L_r(\phi)^{-1}, \quad \mbox{ for all } x \in \mathcal{S}_{\zeta}(\mathbf{\Sigma}).$$
The assignation $\phi \mapsto L_r(\phi)$ defines a projective representation $L_r: G \to \PGL(V)$ and we thus have a procedure $r \mapsto L_r$ to associate to certain representations of stated skein algebras some projective representations of (subgroups of) the mapping class group. The same procedure works for the reduced stated skein algebra as well. 

\underline{\textbf{Irreducible representations fixed by $G$:}}
Theorem \ref{theorem_AL} produces a lot of invariant representations: suppose that $\widehat{x} \in \widehat{X}(\mathbf{\Sigma})$ is such that $(1)$ $\widehat{x}\in \mathcal{AL}(\mathbf{\Sigma})$ and $(2)$ $\widehat{x}$ is fixed by $G$, then the Azumaya representation $r_{\widehat{x}}: \mathcal{S}_{\zeta}(\mathbf{\Sigma})\to \End(V)$ with full shadow $\widehat{x}$ is fixed by $G$ and thus defines a representation $L_{\widehat{x}}: G\to \PGL(V)$. Aside from the whole mapping class group, we will consider two subgroups $G\subset \Mod(\Sigma)$. If $\Sigma$ possess $n\geq 2$ boundary components, a mapping class permutes these boundary components thus defines a permutation in $\mathbb{S}_n$: we denote by $\Mod(\Sigma)^{\partial}$ the kernel of the group morphism $\Mod(\Sigma)\to \mathbb{S}_n$. A mapping class also acts on the homology $\mathrm{H}_1(\Sigma; \mathbb{Z})$ and preserves the intersection form $(\cdot, \cdot)$ so it defines a symplectomorphism; the \textbf{Torelli subgroup} $\mathcal{T}(\Sigma)$  is the kernel of the morphism $\Mod(\Sigma) \to \mathrm{Sp}(\mathrm{H}_1(\Sigma; \mathbb{Z}), (\cdot, \cdot) )$. 

 Here are some examples: 
\begin{itemize}
\item The \textbf{trivial representation} $\rho_0: \pi_1(\Sigma, \mathbb{V}) \to \SL_2$ is the representation sending every path to the identity matrix. It is clearly fixed by the mapping class group so if $\widehat{x}\in \mathcal{AL}(\mathbf{\Sigma})$ (or $\widehat{x}\in \overline{\mathcal{AL}}(\mathbf{\Sigma})$)  is a lift of $\rho_0$, then it will be fixed by $\Mod(\Sigma)^{\partial}$ and thus we obtain a representation $L_{\widehat{x}}: \Mod(\Sigma)^{\partial} \to \PGL(V_{\widehat{x}})$. If moreover all puncture invariants $h_{p_i}$ in $\widehat{x}$ are pairwise equal, we even get a representation  $L_{\widehat{x}}: \Mod(\Sigma) \to \PGL(V_{\widehat{x}})$. By Theorem \ref{theorem_AL} such representations exist for stated skein algebras of the marked surface $\mathbf{\Sigma}_{g,n}^*$: they induce the Lyubashenko representations. Such representations can also be defined using the reduced stated skein algebras of essential marked surfaces $\mathbf{\Sigma}$ which either have a boundary component with at least two boundary edges or which do not have any inner puncture: in this case we can chose $\widehat{x}=(\rho_0, -2, \ldots, -2, h_{\partial_j})$ with all puncture invariants $h_{p_i}=-2$ and the $h_{\partial_j}$ arbitrary. The braid group representations underlying Kashaev link invariant are particular cases of such representations.
\item If $\rho_D: \pi_1(\Sigma, \mathbb{V}) \to \SL_2$ is a diagonal representation (its image lies in the subgroup of diagonal matrices), then $\rho_D$ is fixed by the Torelli group $\mathcal{T}(\Sigma)$. So any element $\widehat{x}$ which is a lift of a diagonal representation $\rho_D$ and which lies in the Azumaya locus induces a representation $L_{\widehat{x}}: \mathcal{T}(\Sigma)\to \PGL(V_{\widehat{x}})$. When $\mathbf{\Sigma}$ is unmarked, such representations were considered in \cite{BCGPTQFT}. 
\end{itemize}

%\underline{\textbf{Explicit constructions of Azumaya representations:}}
%When $\widehat{x}\in \mathcal{AL}(\mathbf{\Sigma})$, Theorem \ref{theorem_AL} proves the existence and unicity up to isomorphism of an Azumaya representation $r_{\widehat{x}}: \mathcal{S}_{\zeta}(\mathbf{\Sigma})\to \End(V_{\widehat{x}})$ with full shadow $\widehat{x}$, however it does not provide an explicit construction of $V_{\widehat{x}}$. Such constructions can be realized using stated skein modules of handlebodies. First suppose that $\mathbf{\Sigma}=(\Sigma_{g},\emptyset)$ is closed, unmarked of genus $g$ and let $H_g$ be a genus $g$ handlebody bounded by $\Sigma_{g}$. By gluing $\Sigma\times [0,1]$ to $H_g$ and retracting, we obtain another copy of $H_g$; this endows $\mathcal{S}_{\zeta}(H_g)$ with a structure of (left) $\mathcal{S}_{\zeta}(\Sigma_g)$ module. Let $\rho: \pi_1(\Sigma_g)\to \SL_2$ be a non central representation and consider the composition $\partial \rho: \pi_1(H_g) \to \pi_1(\Sigma_g) \xrightarrow{\rho} \SL_2$
% with corresponding maximal ideal $m_{[\partial \rho]}\subset \mathcal{S}_{+1}(H_g)$ and consider the quotient: 
%$$ \mathcal{S}_{\zeta}(H_g)_{| [\partial \rho]} := \quotient{\mathcal{S}_{\zeta}(H_g)}{Fr (m_{[\rho]}) \mathcal{S}_{\zeta}(H_g)}.$$
%It is proved in \cite{FrohmanKaniaLe_UnicityRep, KojuKaruo_Azumaya} that $\dim\left( \mathcal{S}_{\zeta}(H_g)_{| [\rho]} \right)$ is the PI-degree of $\mathcal{S}_{\zeta}(\Sigma_g)$, therefore 

 \subsection{Witten-Reshetikhin-Turaev representations}
 
 When $\mathbf{\Sigma}=(\Sigma_g, \emptyset)$ is a closed unmarked surface, the trivial representation $[\rho_0]$ does not belong to the Azumaya locus. Still we can construct an irreducible representation with classical shadow $[\rho_0]$ which is fixed by $\Mod(\Sigma_g)$ as follows. Consider a Heegaard splitting of the sphere $S^3=H'_g \cup_{\Sigma_g} H_g$;  so $H_g, H_g' \subset S^3$ are two genus $g$ handlebodies which intersect along their boundary $H'_g \cap H_g= \partial H_g =\partial H'_g=\Sigma_g$. Denote by $\iota: H_g\hookrightarrow S^3$ and $\iota': H_g'\hookrightarrow S^3$ the inclusion maps and define a pairing via:
 $$ \left(\cdot, \cdot \right)^H: \mathcal{S}_{\zeta}(H_g') \otimes \mathcal{S}_{\zeta}(H_g) \to \mathcal{S}_{\zeta}(S^3) \cong \mathbb{C}, \quad (T', T)^H:= \left< \iota'(T')\cup \iota(T) \right>.$$
 The WRT module is then defined via: 
 $$ V_{\zeta}^{WRT}(\Sigma_g):= \quotient{ \mathcal{S}_{\zeta}(H_g)}{ \mathrm{RKer}( (\cdot, \cdot)^H)}.$$
 The natural action of $\mathcal{S}_{\zeta}(\Sigma_g)$ on $\mathcal{S}_{\zeta}(H_g)$ passes to the quotient and equipped $V^{WRT}_{\zeta}(\Sigma_g)$ with a structure of left $\mathcal{S}_{\zeta}(\Sigma_g)$ module referred as the \textbf{Witten-Reshetikhin-Turaev (WRT) skein representation}. The WRT skein representation is irreducible with classical shadow the trivial representation (\cite{BonahonWong4}). It is a non trivial fact, derived from TQFT (see \cite{RT, Tu,  BHMV2}), that the WRT skein representation is fixed by $\Mod(\Sigma_g)$ and therefore defines a representation 
 $$ L_{WRT}: \Mod(\Sigma_g) \to \PGL(V_{\zeta}^{WRT}(\Sigma_g)).$$
 Such representations also exist for open surfaces where each boundary component is "colored" by some simple $U_q\mathfrak{sl}_2$ modules (\cite{RT, Tu,  BHMV2}).
A lot is known about these representations. Without trying being exhaustive, I   list a few chosen properties and applications. Here we consider $L_{WRT}$ as a linear representation of a central extension $\widetilde{\Mod}(\Sigma_g)$  of $\Mod(\Sigma_g)$.
\par \underline{\textbf{Kernel:}} For $T_{\gamma}\in \Mod(\Sigma_g)$ the Dehn twist around a simple closed curve $\gamma\subset \Sigma_g$, then $L_{WRT}( T_{\gamma}^N)=\id$. Thus $L_{WRT}$ is far from being faithful. Whether the kernel of $L_{WRT}$ only consists in the normal subgroup generated by elements $T_{\gamma}^N$ is an open (and difficult) question; see \cite{DetcherrySantharoubane_KerneWRTRep} for recent progresses. 
\par \underline{\textbf{Infinite image:}} When $g\geq 2$, the image of $L_{WRT}$ is infinite (\cite{Fu99}, see also  \cite{Ma99, Koju3}). In particular the quotient $\quotient{\Mod(\Sigma_g)}{(T_{\gamma}^N, \gamma \subset \Sigma)}$ is infinite as well. It is conjectured that the image of pseudo-Anosov elements have infinite order (see \cite{DetcherryBelettiKaflfagianniYang} for recent progresses). See also \cite{FunarEyssidieux_OrbifoldKahler} for a geometric criterion to prove the infiniteness of the image of a subgroup of $\Mod(\Sigma_g)$. 
\par \underline{\textbf{Decomposition into irreducible:}} $V_{\zeta}^{WRT}(\Sigma_g)$ is semi-simple as a $\widetilde{\Mod}(\Sigma_g)$ module (\cite{BHMV2}). When $g=1$, an explicit decomposition into irreducible factors is well known (\cite{Koju1}). When $N$ is prime, $V_{\zeta}^{WRT}(\Sigma_g)$ is simple (\cite{Ro}). When $4$ divides $N$, $V_{\zeta}^{WRT}(\Sigma_g)$ decomposes as the sum of two submodules (\cite{BHMV2}), they are simple when $g=2$ and $N=4r$ with $r$ an odd prime (\cite{Koju2}). When $g=2$ and $N$ is the product of two (not necessarily distinct) primes  $V_{\zeta}^{WRT}(\Sigma_2)$ is simple (\cite{Koju2}). When $\Sigma_{g,n}$ admits at least one boundary component colored by the standard $2$ dimensional $U_q\mathfrak{sl}_2$ module $S_1$, then $V_{\zeta}^{WRT}(\Sigma_{g,n}, S_1)$ is simple \cite{KoberdaSantharoubane17}. 
\par \underline{\textbf{Asymptotic faithfulness:}} If $\{\zeta_N\}_N$ is an infinite family of roots of unity with unbounded orders, then the direct sum $\oplus_{N} V_{\zeta}^{WRT}(\Sigma_g)$ is a faithful projective module (\cite{And06, MN}).
\par \underline{\textbf{Integrality:}} Let $\mathcal{O}=\mathbb{Z}[\zeta]$ be the cyclotomic ring and suppose that $N$ is prime. Then $V_{\zeta}^{WRT}(\Sigma_g)$ contains a $\mathcal{O}$-lattice $\mathcal{L} \subset V_{\zeta}^{WRT}(\Sigma_g)$ which is preserved by $\Mod(\Sigma_g)$ (\cite{GilmerMasbaum_IntegralLattices}). Let $h:=1-\zeta$ be the (only) maximal ideal of $\mathcal{O}$ and consider the finite quotients $F_{N,k}:= \quotient{\mathcal{O}}{(h^k)}$. Then $\mathcal{L}_{N,k}:=\mathcal{L}\otimes_{\mathcal{O}}F_{N,k}$ is a $\widetilde{\Mod}(\Sigma_g)$ module and thus defines a representation $R_{N,k}: \widetilde{\Mod}(\Sigma_g)\to \GL_{F_{N,k}}(\mathcal{L}_{N,k})$ with finite image. Similar representations also exist for open surfaces.  These finite representations are interesting tools to study $\Mod(\Sigma_g)$; here are two applications.
\begin{itemize}
\item Using the asymptotic faithfulness of the WRT representations, one can show that $\widetilde{\Mod}(\Sigma_g)$ is "well approximated" by its finite quotients $\quotient{\widetilde{\Mod}(\Sigma_g)}{\ker(R_{N,k})}$: this provides an original proof of the fact that $\Mod(\Sigma_g)$ is residually finite (\cite{Funar_MCG_ResFinite}). 
\item When $\chi(\Sigma)<0$, using the Birman exact sequence $0\to \pi_{1}(\Sigma) \to \Mod^1(\Sigma) \to \Mod(\Sigma) \to 0$ (here $\Mod^1(\Sigma)$ is the subgroup of mapping classes fixing some framed base point $v\in \Sigma$) and finite image WRT type representations $\Mod^1(\Sigma)\to \PGL(\mathcal{L}_{N,k})$, we define finite index subgroups $K_{N,k}:= \ker \left(\pi_1(\Sigma_g) \to \Mod^1(\Sigma)\to \PGL(\mathcal{L}_{N,k}) \right)$ and thus finite coverings $\widehat{\Sigma}_{N,k} \to \Sigma$. Such coverings were used in \cite{KoberdaSantharoubane_FiniteCoversHomology} to answer the following question raised by Goldman. Given  $\widehat{\Sigma}\to \Sigma$ a finite (regular) covering of $\Sigma$, a simple closed curve $\gamma \subset \Sigma$ induces an element $\gamma \in \pi_1(\Sigma)$ such that  for $n$ big enough $\gamma^n$ lifts to $\Sigma'$ and thus defines a class $[\gamma^n] \in \mathrm{H}_1(\Sigma'; \mathbb{Z})$.  The question is whether such classes $[\gamma^n]$ generate the whole $\mathrm{H}_1(\Sigma'; \mathbb{Z})$ or not. The authors of \cite{KoberdaSantharoubane_FiniteCoversHomology} proved that this is not case for some coverings  $\widehat{\Sigma}_{N,k} \to \Sigma$.
 \end{itemize}
 
 \subsection{Blanchet-Costantino-Geer-Patureau Mirand representations}
 
 Let $\Sigma=\Sigma_g$ be closed and 
recall the inclusion $\mathrm{H}^1(\Sigma; \mathbb{C}^*) \to \mathcal{X}(\Sigma), \omega \mapsto [\rho_{\omega}]$ whose image consists in diagonal representations. By Theorem \ref{theorem_AL}, for $\omega \in \mathrm{H}^1(\Sigma; \mathbb{C}^*)\setminus \mathrm{H}^1(\Sigma; \mathbb{Z}/2\mathbb{Z})$, the class $[\rho_{\omega}]$ belongs to the Azumaya locus of $\mathcal{S}_{\zeta}(\Sigma)$;  we therefore have an Azumaya representation $r_{\omega}: \mathcal{S}_{\zeta}(\Sigma) \to \End(V^{BCGP}(\Sigma, \omega))$ with full shadow $[\rho_{\omega}]$. By definition, elements of the Torelli group act trivially on homology so $r_{\omega}$ is fixed by $\mathcal{T}(\Sigma)$ and $r_{\omega}$ defines a projective representation
$$ L_{r_{\omega}}: \mathcal{T}(\Sigma) \to \PGL(V^{BCGP}(\Sigma, \omega)).$$
These representations were first defined in \cite{BCGPTQFT} where the authors proved that the image by $L_{r_{\omega}}$ of Dehn twists have infinite order (unlike the WRT representations) and conjectured that $L_{r_{\omega}}$ is injective. This conjecture would imply the projective linearity of $\mathcal{T}(\Sigma)$. Despite of the importance of such a consequence, as far as I know these representations have been poorly studied and no significative progresses have been made towards this conjecture. 
 
 \subsection{Braid groups representations underlying quantum invariants of links}
 For $n\geq 0$, the $n$-th punctured bigon $\mathbb{D}_n=(D_n, \{b_L, b_R\})=\adjustbox{valign=c}{\includegraphics[width=0.7cm]{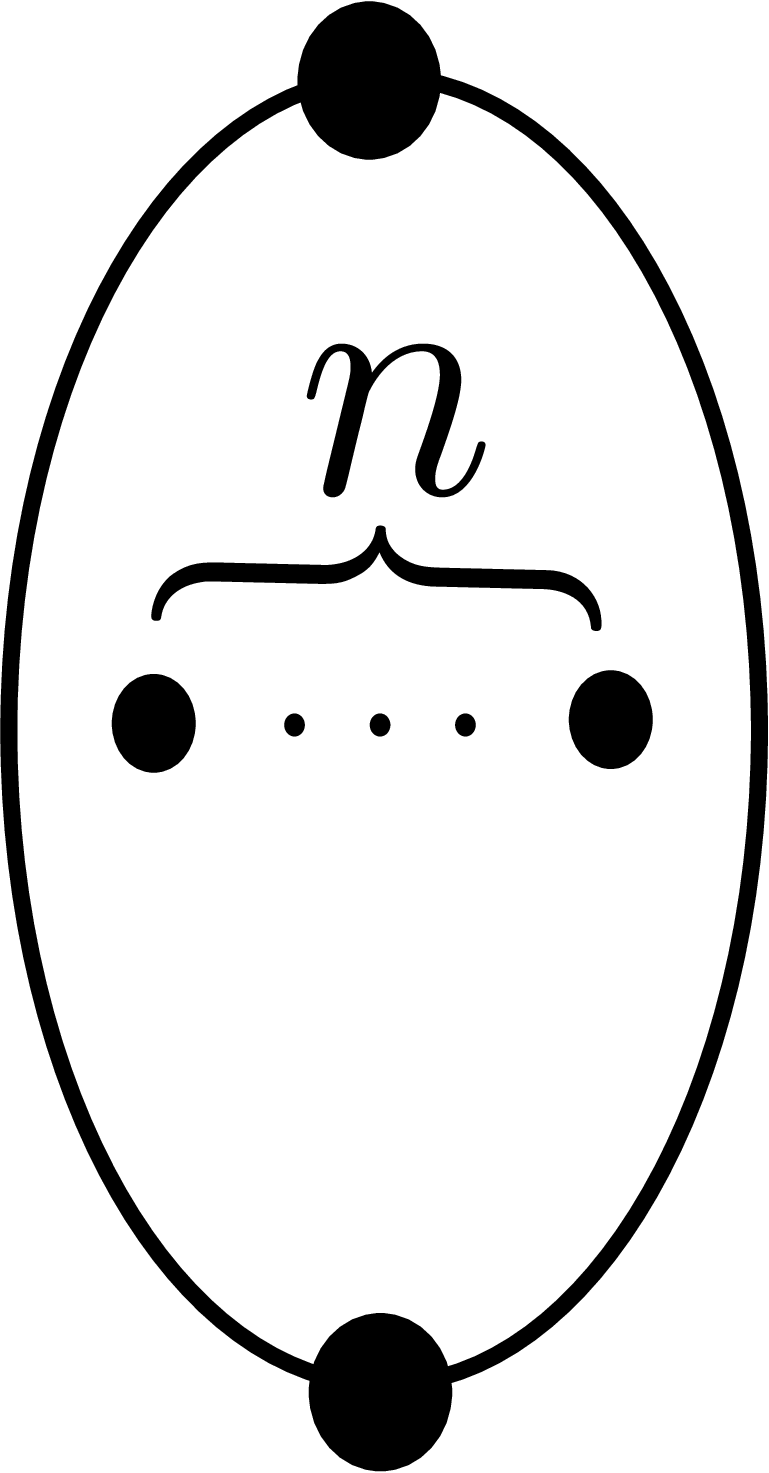}} $ is a disc with two boundary edges and $n$ open subdiscs removed from its interior (so it has two boundary punctures and $n$ inner punctures). Let $\widehat{x}=(\rho, h_{p_1}, \ldots, h_{p_n}, h_{\partial}) \in \widehat{\overline{\mathcal{X}}}(\mathbb{D}_n)$. By Theorem \ref{theorem_AL}, $\widehat{x} \in \overline{\mathcal{AL}}(\mathbb{D}_n)$ if and only if for all $1\leq i \leq n$ then either $\rho(\gamma_{p_i})\neq \pm \mathds{1}_2$ or $h_{p_i}=\pm 2$; in this case we have an Azumaya irreducible representation $r_{\widehat{x}}: \overline{\mathcal{S}}_{\zeta}(\mathbb{D}_n)\to \End(V_{\widehat{x}})$ with full shadow $\widehat{x}$. 
 \par Let us first assume that $\widehat{x}_0=(\rho_0, -2, \ldots, -2, 1)$ is the trivial representations with punctures invariants $h_{p_i}=-2$ and boundary invariant $h_{\partial}=1$. In this case we obtain a representation of the braid group:
 $$ L_{r_{\widehat{x}_0}}: B_n \to \PGL(V_{\widehat{x}_0})$$
 which $(1)$ can be linearized and $(2)$ can be extended to a linear representation $L_0: fB_n \to \GL(V_{\widehat{x}_0})$ of the framed braid group $fB_n$ by means of twist operators.  For $L \subset S^3$ a framed link which is the Markov closure of a framed braid $\beta \in fB_n$, the \textbf{Kashaev invariant} $\left< L \right>$  is recovered as  the trace
 $$ \left< L \right>:= \tau^R ( L_0(\beta))$$
 where $\tau^R$ is the so-called \textbf{renormalized trace} (see below).
 \par  
 Now let $\beta \in fB_n$ be a framed braid which fixes $\widehat{x}$ (i.e. $\beta_* \widehat{x}=\widehat{x}$). We can then define an intertwiner $L_{r_{\widehat{x}}}(\beta) \in \GL(V_{\widehat{x}})$, uniquely defined up to multiplication by a $N^2$ root of unity, and 
 when $\widehat{x}=(\rho, h_{p_1}, \ldots, h_{p_n}, h_{\partial})\in \overline{\mathcal{AL}}(\mathbb{D}_n)$ is such that none of the $\rho(\gamma_{p_i})$ is parabolic, we can 
 consider the trace $\tau^R(L_{r_{\widehat{x}}}(\beta) )$. This trace defines an invariant of triples $\mathbb{L}=(L, [\rho_L], h)$ where: 
 \begin{itemize}
 \item
  $L=L_1 \sqcup \ldots \sqcup L_k\subset S^3$ is the framed link which is the Markov closure of $\beta$; we denote by $M_L:=S^3 \setminus N(L)$ its exterior. 
  \item $[\rho_L]$ is the conjugacy class of the representation $\rho_L: \pi_1(M_L) \to \SL_2$ induced from the classical shadow $\rho: \pi_1(\mathbb{D}_n, \mathbb{V})\to \SL_2$ of $V_{\widehat{x}}$ by the natural morphism $\pi: \mathbb{D}_n \to M_L$.
  \item  $h=(h_{L_1}, \ldots, h_{L_n})$ corresponds to the puncture invariants $h_{p_i}$ under the natural map induced by $\pi$ which associates to each inner puncture $p$ its corresponding component $L_p \subset L$.
  \end{itemize}
  The invariant $\left< \mathbb{L} \right>$ is well defined up to multiplication by a $N^2$-th root of unity (i.e. does not depend on the choices of $\beta$, $\rho$ and $h_{\partial}$) and recover many "quantum groups type invariants": 
 \begin{enumerate}
 \item As we just said, when $[\rho_0]$ is the class of the trivial representation and $h_{L_i}=-2$ for all $i$, $\left< L, [\rho_0], h\right>= \left< L \right>$ is the Kashaev invariant.
 \item When $\rho_{\omega} : \pi_1(M_L) \to \SL_2$ is a diagonal representation induced by a cohomology class $\omega \in \mathrm{H}^1(M_L; \mathbb{C}^*)$, then $\left< L, [\rho_{\omega}] \right>$ is the link invariant defined in \cite{GeerPatureau_LinksInv}.
 \item The link invariant $\left< \mathbb{L} \right>$  for non diagonal $\rho$ where first defined in \cite{BGPR_Biquandle}.
 \end{enumerate}
 For details concerning this skein approach to link invariants see \cite{KojuQGroupsBraidings}; let me just say a few words about the relation with quantum groups. Consider the strong embedding $\iota: \mathbb{D}_1\to \mathbb{D}_n$ which consists in placing the $n$ inner punctures of $\mathbb{D}_n$ inside the only inner puncture of $\mathbb{D}_1$; so we get an embedding $\iota_*: \overline{\mathcal{S}}_{\zeta}(\mathbb{D}_1) \to \overline{\mathcal{S}}_{\zeta}(\mathbb{D}_n)$. By gluing $n$ copies of $\mathbb{D}_1$ together we obtain $\mathbb{D}_n$ so we have a splitting morphism $\theta_n: \overline{\mathcal{S}}_{\zeta}(\mathbb{D}_n) \to \overline{\mathcal{S}}_{\zeta}(\mathbb{D}_1)^{\otimes n}$. The composition
 $$ \Delta: \overline{\mathcal{S}}_{\zeta}(\mathbb{D}_1) \xrightarrow{\iota_*} \overline{\mathcal{S}}_{\zeta}(\mathbb{D}_2) \xrightarrow{\theta_2} \overline{\mathcal{S}}_{\zeta}(\mathbb{D}_1)^{\otimes 2}$$
 defines a coproduct which endows $\overline{\mathcal{S}}_{\zeta}(\mathbb{D}_1)$ with a structure of Hopf algebra isomorphic to the Drinfeld double $D_qB$ of the quantum Borel algebra $B_q$. In particular the quotient $\quotient{\overline{\mathcal{S}}_{\zeta}(\mathbb{D}_1)}{(\alpha_{\partial}-1)}$ is isomorphic to the simply connected quantum group $U_q\mathfrak{sl}_2$. The category $\mathcal{C}$ of $D_qB$-modules on which $Z^0_{\mathbb{D}_1}$ acts semi-simply admits a pivotal structure and its subcategory $\Proj(\mathcal{C})$ of projective objects admits a trace $\tau^R : \mathrm{HH}_0(\Proj(\mathcal{C}))\to \mathbb{C}$ studied in \cite{GeerPatureau_TraceQG}. If $V_1, \ldots, V_n$ are simple $D_qB$ modules, then using $\theta_n$, the tensor product $V_1\otimes \ldots \otimes V_n$ admits a structure of simple $\overline{\mathcal{S}}_{\zeta}(\mathbb{D}_n)$ module. Conversely, it is proved in \cite{KojuKaruo_RepRSSkein} that every simple $\overline{\mathcal{S}}_{\zeta}(\mathbb{D}_n)$ module arises that way; so every Azumaya representations can be written as $V_{\widehat{x}}= V_1\otimes \ldots \otimes V_n$. Using $\iota_*$, such an Azumaya representation $V_{\widehat{x}}$ acquires a structure of $D_qB$ module and the fact that none of the $\rho(\gamma_{p_i})$ is parabolic implies that $V_{\widehat{x}}\in \Proj(\mathcal{C})$. By design, the action of the braid group on $\mathbb{D}_n$ leaves the image of $\iota$ invariant so $\beta_* \iota_*(x)= \iota_*(x)$ for all $x\in D_qB$. This implies that  $L_{r_{\widehat{x}}}(\beta): V_{\widehat{x}}\to V_{\widehat{x}}$ is an equivariant morphism of $D_qB$ modules so an endomorphism of $\Proj(\mathcal{C})$ to which we can apply $\tau^R$. 
 \par Note that the pure braid groups representations underlying the colored Jones polynomials admit a similar skein description. Let $S_i$ be the $i+1$-th dimension simple $U_q\mathfrak{sl}_2$ module ($i<N-1$) and consider the $\overline{\mathcal{S}}_{\zeta}(\mathbb{D}_n)$ simple module $V= S_{i_1}\otimes \ldots \otimes S_{i_n}$ for some $1\leq i_k \leq N-1$. Its classical shadow is the trivial representation $\rho_0$ and its full shadow is not in the Azumaya locus. Still $V$ is stable under the action of the pure braid group $PB_n \subset B_n$ thus defines a projective representation $L_V: PB_n \to \PGL(V)$ which can be linearized and extended to representation $L_V: fPB_n \to \GL(V)$ of the framed pure braid group. This time an intertwiner $L_V(\beta)$ defines an endomorphism of $\mathcal{C}$, but not of $\Proj(\mathcal{C})$, and we need to consider a different trace on $\mathcal{C}$ named the \textbf{q-trace} $\qtr$. 
  If $L$ is the Markov closure of $\beta$ and $\beta(i_k)=i_k$ for all $k$, then the colored Jones polynomial is the trace $J (L) = \qtr(L_V(\beta))$. 
 
 \subsection{Lyubashenko representations}
 
 By Theorem \ref{theorem_AL} the trivial representation $\rho_0: \pi_1(\Sigma_{g,1}) \to \SL_2$ (sending every paths to the identity matrix) belongs to the Azumaya locus of $\mathcal{S}_{\zeta}(\mathbf{\Sigma}_{g, 0}^*)$, therefore we have a projective representation
 $$ L_0: \Mod(\Sigma_{g,1}) \to \PGL(V_{\rho_0}).$$
 
 \begin{theorem}(\cite{Faitg_LGFT_MCG, Faitg_LGFT_SSkein})\label{theorem_Lyubashenko} $L_0$ is the Lyubashenko representation from \cite{Lyubashenko_InvariantsMCGRep} associated to the small quantum group $\mathfrak{u}_q\mathfrak{sl}_2$. \end{theorem}
 
 The theorem is an immediate consequence of the work in \cite{Faitg_LGFT_MCG, Faitg_LGFT_SSkein} though Matthieu did not formulated the theorem that way. I will first explain why this theorem is merely a reformulation of the main theorem in  \cite{Faitg_LGFT_MCG}. I will explain why this theorem is also a consequence of an interesting relationship found in  \cite{KojuMurakami_QCharVar} between the Kerler-Lyubashenko TQFT of \cite{KerlerLyubashenko_TQFT} and stated skein modules of $3$-manifolds which we now introduce.
 
 \subsubsection{Marked $3$-manifolds} Consider the disc $\mathbb{D}:=\{ (x,y) \in \mathbb{R}^2: x^2+y^2\leq 1\}$. The \textbf{height} of $(x,y)\in \mathbb{D}$ is $h((x,y)):=y$. 
 A \textbf{marked 3 manifold} is a pair $\mathbf{M}=(M, \mathcal{N})$ where $M$ is a compact oriented $3$ manifold and $\mathcal{N}$ is a finite collection of pairwise disjoint embedded discs $\mathbb{D}_i \subset \partial M$ named \textbf{boundary discs}. Each boundary disc is considered parametrized by an oriented homeomorphism $\mathbb{D}\cong \mathbb{D}_i$ so we can define the height of an element in $\mathbb{D}_i$. An \textbf{embedding} $f: (M, \mathcal{N}) \to (M', \mathcal{N}')$ is an oriented embedding $f: M \hookrightarrow M'$ which induces an embedding $f: \mathcal{N} \to \mathcal{N}'$ which is height increasing in the sense that if $v,w \in \mathbb{D}_i$ are such that $h(v)<h(w)$ in $\mathbb{D}_i$ and $f$ sends $\mathbb{D}$ to $\mathbb{D}'$ then $h(f(v))<h(f(w))$ in $\mathbb{D}'$. Also if both $\mathbb{D}_i, \mathbb{D}_j$ are sent into $\mathbb{D}'$, we impose that their image have pairwise disjoint height: $h(f(\mathbb{D}_i))\cap h(f(\mathbb{D}_j))= \emptyset$. We get a category $\mathcal{M}$ of marked $3$-manifolds. A marked surface $(\Sigma, \mathcal{A})$ defines a marked $3$-manifold $(\Sigma\times [0,1], \mathcal{A}\times[0,1])$ so we have the inclusions $\mathrm{MS}^{str}\subset \mathrm{MS} \subset \mathcal{M}$. Stated skein modules $\mathcal{S}_A(M)$ of marked $3$ manifolds are defined in the same manner than for marked surfaces and we get a functor $\mathcal{S}_A: \mathcal{M} \to \Mod_k$. For $a,b$ two boundary discs of $\mathbf{M}$, the marked $3$-manifold $\mathbf{M}_{a \# b}$ is obtained from $\mathbf{M}$ by gluing $a$ and $b$ together and the definition of the splitting morphism $\theta_{a\# b} : \mathcal{S}_A(\mathbf{M}_{a\# b}) \to \mathcal{S}_A(\mathbf{M})$ extends straightforwardly to marked $3$-manifolds (\cite{BloomquistLe20}) with one difference: $\theta_{a\#b}$ may not be injective anymore (\cite{CostantinoLe_SSkeinTQFT}). 
 In particular, $\mathcal{S}_A(\mathbf{M})$ is a $(\mathcal{O}_q\SL_2)^{\otimes \mathcal{N}}$ comodule. Theorem \ref{theorem_fusion} still holds for marked $3$-manifolds (\cite{CostantinoLe_SSkeinTQFT}).
 Also the Frobenius morphism $Fr_{\mathbf{M}}: \mathcal{S}_{+1}(\mathbf{M}) \to \mathcal{S}_{\zeta}(\mathbf{M})$ is still well-defined \cite{BloomquistLe20} (though not necessarily injective \cite{CostantinoLe_SSkeinTQFT}) and its image consists in transparent elements in the sense that for any stated tangles $(T_1,s_1), (T_2,s_2)$ then the class $Fr_{\mathbf{M}}([T_1,s_1])\cup [T_2,s_2] \in \mathcal{S}_{\zeta}(\mathbf{M})$ does not depend on how $T_1$ and $T_2$ are entangled so $Fr$ endows $\mathcal{S}_{\zeta}(\mathbf{M})$ with a structure of $\mathcal{S}_{+1}(\mathbf{M})$ module. In particular for $\mathfrak{m}\in \Specm(\mathcal{S}_{+1}(\mathbf{M}))$ we can consider the quotient
 $$ \mathcal{S}_{\zeta}(\mathbf{M})_{\mathfrak{m}}:= \quotient{\mathcal{S}_{\zeta}(\mathbf{M})}{(Fr_{\mathbf{M}}(\mathfrak{m}) \mathcal{S}_{\zeta}(\mathbf{M}))}.$$
 We will be interested in the subcategory $\mathcal{M}_c^{(1)} \subset \mathcal{M}$ of marked $3$-manifolds $\mathbf{M}=(M, \mathcal{N})$ such that $M$ is connected and $\mathcal{N}$ has a single boundary disc. $\mathcal{M}_c^{(1)}$ has a natural monoidal structure given by the fusion operation: for $\mathbf{M}_1=(M_1,\{a_1\}), \mathbf{M}_2=(M_2, \{a_2\})\in \mathcal{M}_c^{(1)}$ we write $\mathbf{M}_1\wedge \mathbf{M}_2:= (\mathbf{M}_1\sqcup \mathbf{M}_2)_{a_1\circledast a_2} \in \mathcal{M}_{c}^{(1)}$. 
 \par  
 For $\mathbf{M}\in \mathcal{M}_c^1$, then $\mathcal{S}_{+1}(\mathbf{M}) \cong \mathcal{O}[\mathcal{R}_{\SL_2}(M)]$ is isomorphic to the ring of regular function of the representation variety so we can identify a maximal ideal $\mathfrak{m}\in \Specm(\mathcal{S}_{+1}(\mathbf{M}))$ with a representation $\rho: \pi_1(M)\to \SL_2$. When $\rho=\rho_0$ is the trivial representation we will be interested in the quotient $\mathcal{S}^{[0]}_{\zeta}(\mathbf{M}):= \mathcal{S}_{\zeta}(\mathbf{M})_{\rho_0}$.  Given an embedding $f: \mathbf{M}\to \mathbf{M}'$, the morphism $f_*$ passes to the quotient and induces a morphism $f_*: \mathcal{S}_{\zeta}^{[0]}(\mathbf{M})\to \mathcal{S}_{\zeta}^{[0]}(\mathbf{M}')$ so we have a functor $\mathcal{S}^{[0]}_{\zeta}: \mathcal{M}\to \Vect$. 
 For the bigon, the algebra 
 $$ \mathcal{S}^{[0]}_{\zeta}(\mathbb{B}) \cong \quotient{\mathcal{O}_{\zeta}(\mathbb{B})}{(a^N=d^N=1, b^N=c^N=0)} =: \mathfrak{o}_q\SL_2, \quad q:=\zeta^2, $$
 acquires a Hopf algebra structure by passing to the quotient and, as a Hopf algebra, $\mathfrak{o}_q\SL_2$ is isomorphic to the dual of the small quantum group $\mathfrak{u}_q\mathfrak{sl}_2:= \quotient{U_q\mathfrak{sl}_2}{(E^N=F^N=0, K^N=1)}$. Let $\mathcal{C}$ be the braided category of locally free $\mathfrak{u}_q\mathfrak{sl}_2$-modules, or equivalently of $\mathfrak{o}_q\SL_2$- comodules. For $\mathbf{M}\in \mathcal{M}_c^{(1)}$, $\mathcal{S}^{[0]}_{\zeta}(\mathbf{M})$ is a $\mathfrak{o}_q\mathfrak{sl}_2$ comodule, so an element of $\mathcal{C}$. The $3$-manifolds version of Theorem \ref{theorem_fusion} proved in \cite{CostantinoLe_SSkeinTQFT} implies that the functor
 $$\mathcal{S}^{[0]}_{\zeta}: (\mathcal{M}_c^{(1)}, \wedge) \to (\mathcal{C}, \overline{\otimes})$$
 is Lax monoidal (see \cite{KojuMurakami_QCharVar} for details). 
 
 \subsubsection{Proof of Theorem \ref{theorem_Lyubashenko} using Faitg's work} 

Write $\mathcal{H}:= \mathfrak{u}_q\mathfrak{sl}_2$.
 The module $V_{\rho_0}$ can be constructed explicitly as follows. Let $H_g^* \in \mathcal{M}_c^{(1)}$ be a genus $g$ handlebody with one boundary disc (so $H_g^*\cong \mathbf{\Sigma}_{0, g}^* \times [0,1]$). Since $H_g^*$ is a module object over $\mathbf{\Sigma}_{g,0}^*$ in $ \mathcal{M}_c^{(1)}$, by functoriality $\mathcal{S}_{\zeta}^{[0]}(H_g^*)$ is a module over $\mathcal{S}_{\zeta}^{[0]}(\mathbf{\Sigma}_{g,0}^*)$. Using the Daisy graph $\Gamma$ presenting $\mathbf{\Sigma}_{0,g}^*$, we obtain the linear isomorphism 
 $\Psi^{\Gamma}: \mathcal{S}_{\zeta}(H_g^*) \cong (\mathcal{O}_q\SL_2)^{\otimes g}$ which induces the isomorphism $\Psi^{\Gamma}_{[0]}: \mathcal{S}_{\zeta}^{[0]}(H_g^*) \cong (\mathfrak{o}_q\SL_2)^{\otimes g}\cong (\mathcal{H}^*)^{\otimes g}$. In particular $\dim( \mathcal{S}_{\zeta}^{[0]}(H_g^*)) = N^{3g} = PI-deg(\mathcal{S}_{\zeta}(\mathbf{\Sigma}_{g,0}^*))$. Therefore $\mathcal{S}_{\zeta}(H_g^*)$ is an Azumaya representation of $\mathcal{S}_{\zeta}^{[0]}(\mathbf{\Sigma}_{g,0}^*)$ with shadow $\rho_0$ so we can choose $V_{\rho_0}:= \mathcal{S}_{\zeta}^{[0]}(\mathbf{\Sigma}_{g,0}^*)$. 
 
 \vspace{2mm}
 \par In \cite{Faitg_LGFT_MCG, Faitg_LGFT_SSkein}, Faitg considered an algebra $\mathcal{L}_{g,0}(\mathcal{H})$ and equipped $(\mathcal{H}^*)^{\otimes g}$ with a structure of left module $r: \mathcal{L}_{g,0}(\mathcal{H})\to \End((\mathcal{H}^*)^{\otimes g})$. In \cite{Faitg_LGFT_SSkein} he constructed an algebra isomorphism $\mathcal{L}_{g,0}(\mathcal{H}) \cong \mathcal{S}_{\zeta}^{[0]}(\mathbf{\Sigma}_{g,0}^*)$ (both are isomorphic to $\Mat_{N^{3g}}(\mathbb{C})$ but the isomorphic is explicit) in such a way that the linear isomorphism $\Psi^{\Gamma}_{[0]}: \mathcal{S}_{\zeta}^{[0]}(H_g^*) \cong  (\mathcal{H}^*)^{\otimes g}$ is equivariant. In \cite{Faitg_LGFT_MCG} he proved that the Lyubashenko representation $L^{Lyub}: \Mod(\Sigma_{g,1}) \to \PGL( (\mathcal{H}^*)^{\otimes g})$ is characterized as being the intertwiner 
 $$ r(\phi^{-1}(x)) = L^{Lyub}(\phi) r(x) L^{Lyub}(\phi)^{-1} , \quad \mbox{ for all }\phi \in \Mod(\Sigma_{g,1}), x \in \mathcal{L}_{g,0}(\mathcal{H}).$$
 Theorem \ref{theorem_Lyubashenko} is an immediate consequence of these two theorems.
 
 \subsubsection{Kerler-Lyubashenko TQFT and stated skein module}.
 \par 
 \underline{\textbf{Cobordism category:}}
 The \textbf{Crane-Yetter-Kerler cobordism category} $\Cob$ is the category with objects the surfaces $\Sigma_{g,1}$ for $g\geq 0$ and morphisms the connected compact $3$-manifolds with corners seen as cobordisms. A monoidal structure is given by gluing the boundaries of two surfaces $\Sigma_{g_1, 1}$ and $\Sigma_{g_2,1}$ to two boundary components of a pair of pants to get the surface $\Sigma_{g_1,1}\wedge \Sigma_{g_2,1}=\Sigma_{g_1+g_2, 1}$; so $\Cob$ is a PROP, i.e. a category with set of objects $\mathbb{N}$ and tensor product $g_1 \wedge g_2=g_1+g_2$. $\Cob$ has a natural braided balanced structure coming the natural action of the framed little disc operad where the action of a configuration of $n$ enumerated small  discs inside a big one on $n$ surfaces $\Sigma_{g_1, 1}, \ldots, \Sigma_{g_n,1}$ is defined by gluing each $\Sigma_{g_i, 1}$ to the $i$-th small disc to get a surface $\Sigma_{g_1+\ldots +g_n, 1}$. 
 There is a monoidal functor (actually braided) 
 $$ \partial : \mathcal{M}_c^{(1)} \to \Cob$$
 which sends $\mathbf{M}=(M, \{\mathbb{D}_M\})$ to $\partial \mathbf{M}:= \partial M \setminus \mathring{\mathbb{D}}_M \in \Cob$ and sends an embedding $f: \mathbf{M} \to \mathbf{M}'$ to  $\partial f : \partial \mathbf{M} \to \partial \mathbf{M}'$ where $\partial f= M' \setminus \mathring{f(M)}$. Let $\BT\subset \mathcal{M}_c^{(1)}$ be the full subcategory generated by the objects $H_g^*$ for $g\geq 0$. $\BT$ is a PROP as well and the restriction $\partial: \BT \to \Cob$ is  faithful but not full as we shall see.  
 \par 
 \underline{\textbf{Signature defect:}} $\Cob$ admits a "central extension" $\widetilde{\Cob}$ defined as follows. For each $g$, we have fixed a handlebody $H_g^*=(H_g, \{\mathbb{D}_g\})$ and a boundary $\Sigma_{g,1}=\partial H_g \setminus \mathring{\mathbb{D}_g}$. 
 Write $\Sigma_g:= \partial H_g$ and 
 consider the Lagrangian $\mathcal{L}_g \subset \mathrm{H}_1(\Sigma_g; \mathbb{Q})$ defined by $\mathcal{L}_g:= \ker ( \mathrm{H}_1(\Sigma_g; \mathbb{Q}) \to \mathrm{H}_1(H_g; \mathbb{Q}))$. A cobordism $M \in \Cob(g_1, g_2)$ defines Lagrangians $M_* \mathcal{L}_{g_1} \subset \mathrm{H}_1(\Sigma_{g_2}; \mathbb{Q})$ and $M^* \mathcal{L}_{g_2} \subset \mathrm{H}_1(\Sigma_{g_1}; \mathbb{Q})$ in the obvious way. $\widetilde{\Cob}$ has the same objects than $\Cob$ (i.e. $\mathbb{N}$) with morphisms $\widetilde{\Cob}(g_1, g_2) := \Cob(g_1, g_2) \times \mathbb{Z}$ and composition law for $(M_1, n_1) \in \Cob(g_1, g_2)$ and $(M_2, n_2) \in \Cob(g_2, g_3)$ given by 
 $$ (M_1,n_1) \circ (M_2, n_2) := (M_1\circ M_2, n_1 + n_2 + \mu)$$
 where $\mu = \mu((M_1)_* \mathcal{L}_{g_1}, \mathcal{L}_{g_2}, (M_2)^* \mathcal{L}_{g_3})$ is the Maslov index of the three Lagrangians (see \cite{Tu}). There is a natural forgetful functor $\widetilde{\Cob}\to \Cob$ sending $(M,n)$ to $M$ and the functor $\partial$ factorizes as $\BT \xrightarrow{\widetilde{\partial}} \widetilde{\Cob} \to \Cob$ where $\widetilde{\partial}$ sends $H_g^*$ to $\Sigma_{g,1}$ and an embedding $f: H_{g_1}^* \to H_{g_2}^*$ to $(\partial f, 0)$.

 \par 
 \underline{\textbf{Bottom tangles:}}
 \par Embeddings $H_g^* \to \mathbf{M}$ can be depicted using \textbf{bottom tangles}. A bottom tangle in $\mathbf{M}\in \mathcal{M}_c^{(1)}$ is simply a tangle $T \subset M$ (as we defined before) with the special property that if $T=T_1\sqcup \ldots \sqcup T_g$ are the connected components of $T$ then (up to reindexing) the endpoints of $T_i$ have higher heights (in $\mathbb{D}_M$) than the endpoints of $T_j$ whenever $i<j$. Consider a tubular neighborhood $N(T)$ of $T\cup \mathbb{D}_M$.  Using the framing of $T$, we obtain an identification $H_g \cong N(T)$ so we have an embedding $H_g^* \to \mathbf{M}$ well defined up to isotopy. Conversely taking the spine of the image of an embedding $H_g^* \to \mathbf{M}$ defines a bottom tangle and we have a $1:1$ correspondence between isotopy classes of bottom tangles in $\mathbf{M}$ with $g$ components and $\mathcal{M}(H_g^*, \mathbf{M})$. Let us identify a $g$ bottom tangle $T_g\subset H_{g'}$, its corresponding embedding $i \in \BT(g,g')$ and the cobordism $\widetilde{\partial} i \in \widetilde{\Cob}(g,g')$. 
  \par 
 \underline{\textbf{BPH and BTH algebras:}} What initiated the story of this subsection is the discovery, made independently in \cite{CraneYetter_Categorification} and \cite{Kerler_QInv}, that $\Sigma_{g,1}$ is a braided Hopf algebra object in $\widetilde{\Cob}$. It is even better: it is a Hopf algebra object with cotwist, integral and cointegral. The product $\mu$, unit $\eta$, coproduct $\Delta$, counit $\epsilon$, invertible antipode $S^{\pm 1}$, cotwist $\theta^{\pm 1}$, integral $\lambda$ and cointegral $\Lambda$ are defined by the following bottom tangles: 
 $$ \adjustbox{valign=c}{\includegraphics[width=12cm]{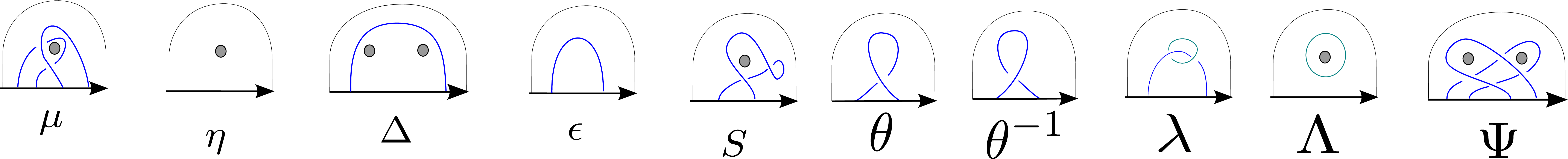}} $$
 where $\Psi$ is the braiding and the green curves in the integral and cointegral mean we should perform surgery along these green curves to obtain the desired cobordisms; so $\lambda, \Lambda$ are morphisms in $\Cob$ but not in $\BT$ whereas all other morphisms belong to $\BT$. In particular $H_1^* \in \BT$ is a Hopf algebra object with cotwist as well but has no integral/cointegral. In \cite{BobtchevaPiergallini}, the authors proved that $\Cob$ is finitely presented by the generators $\{\mu, \eta, \Delta, \epsilon, S^{\pm 1}, \theta^{\pm 1}, \lambda, \Lambda\}$ together with the braiding $\Psi$,  with relations the usual ones defining a (braided) Hopf algebra with cotwist, integral and cointegral plus finite  additional relations. For $\mathcal{C}$ a monoidal category, we will call a \textbf{BPH algebra} the data of a symmetric monoidal functor $Q_H: \Cob \to \mathcal{C}$ or, equivalently, the data $(H, \mu, \eta, \Delta, \epsilon, S^{\pm 1}, \theta^{\pm 1}, \lambda, \Lambda)$ of a Hopf algebra object in $\mathcal{C}$ with cotwist, integral, cointegral and which satisfies the additional BP relations. 
 \par Similarly, it is proved in \cite{KojuMurakami_QCharVar} that $\BT$ is generated by $(H_1^*, \mu, \eta, \Delta, \epsilon, S^{\pm 1}, \theta^{\pm 1})$ though no complete set of relations is known. We define a \textbf{BTH algebra} in $\mathcal{C}$ in the same manner:  as monoidal functors $Q_H: \BT \to \mathcal{C}$, i.e. as Hopf algebra object with cotwists $(\mu, \eta, \Delta, \epsilon, S^{\pm 1}, \theta^{\pm 1})$ satisfying the needed relations imposed by $\BT$. Clearly BPH algebras are BTH algebras, but the converse is not true. For instance, the very existence of the stated skein module functor $\mathcal{S}_A: \BT \to \mathcal{C}=\mathcal{O}_q\SL_2-\Comod$ proves that the braided quantum group $B_q\SL_2$ (transmutation of $\mathcal{O}_q\SL_2$) is a BTH algebra (regardless of the ground ring $k$). However, its is well known that this algebra does not admit any integral nor cointegral so it does not extend to a BPH algebra. 
   \par 
 \underline{\textbf{Kerler-Lyubashenko TQFT:}} Let $\mathcal{C}=\mathfrak{u}_q\mathfrak{sl}_2-\Mod^{lc}$ and consider the transmutation $H:=B\mathfrak{o}_q\SL_2\in \mathcal{C}$ of $\mathfrak{o}_q\SL_2$. Unlike $B_q\SL_2$, $H$ admits an integral and a cointegral (dual to the ones of $\mathfrak{u}_q\mathfrak{sl}_2$) which endow $H$ with a structure of BPH-algebra; we thus have a unique braided functor 
 $$ F^{KL}: \widetilde{\Cob} \to \mathcal{C}$$
 sending $\Sigma_{g,1}$ to $H$ by preserving the BPH structure morphisms: this is the \textbf{Kerler-Lyubashenko TQFT} \cite{KerlerLyubashenko_TQFT}. A mapping class $\phi\in \Mod(\Sigma_{g,1})$ defines a mapping cylinder (pinched along $\partial \Sigma_{g,1}$) $C(\phi): \Sigma_{g,1} \to \Sigma_{g,1}$ seen as a cobordism in $\Cob$. The group $\widetilde{\Mod}(\Sigma_{g,1}) \subset \widetilde{\Cob}(g,g)$ of elements $(C(\phi), n)$ for $\phi \in \Mod(\Sigma_{g,1})$, $n\in \mathbb{Z}$ is a central extension of $\Mod(\Sigma_{g,1})$ and the restriction of the KL functor:
$$ L^{Lyub}: \widetilde{\Mod}(\Sigma_{g,1}) \to \End(F^{KL}(\Sigma_{g,1})), \quad (C(\phi), n) \mapsto F^{KL}((C(\phi), n))$$ 
is the \textbf{Lyubashenko representation} which first appeared in \cite{Lyubashenko_InvariantsMCGRep}. Let us relate $L^{Lyub}$ to the representation $L_{0}$ from skein theory.  In order to alleviate the notations, let us identify $\mathcal{S}_{\zeta}^{[0]}(\mathbf{H}_1^*)$ and $B \mathfrak{o}_q\SL_2$ (they are isomorphic BPH algebras in $\mathcal{C}$ by Theorem \ref{theorem_fusion}); we thus also identify $\mathcal{S}^{[0]}_{\zeta}(\mathbf{H}_g^*)$ and $F^{FL}(\Sigma_{g,1})$. 
 First note that the following diagram commutes: 
$$ \begin{tikzcd}
\BT \ar[rrd, "\restriction{\mathcal{S}_{\zeta}^{[0]}}{\BT}"]  \ar[d, "\widetilde{\partial}"] &{} & {} \\
\widetilde{\Cob} \ar[rr, "F^{KL}"] & {} & \mathcal{C}
 \end{tikzcd}
 $$
 Indeed, both braided functors $\restriction{\mathcal{S}_{\zeta}^{[0]}}{\BT}$ and $F^{KL}\circ \widetilde{\partial}$ send $H_1^*$ to $\mathfrak{o}_q\SL_2$ and agree on  the generators $(\mu, \eta, \Delta, \epsilon, S^{\pm 1}, \theta^{\pm 1})$ of $\BT$ so they must be equal. 
  So $\mathcal{S}_{\zeta}^{[0]}(H_g^*)$ is a left module both for $\mathcal{S}_{\zeta}(\Sigma_g^*)$, via the representation $r_{0}: \mathcal{S}_{\zeta}(\mathbf{\Sigma}_{g,0}^*) \to \End(\mathcal{S}_{\zeta}^{[0]}(H_g^*))$, and for the group $\widetilde{\Mod}(\Sigma_{g,1})$, via the functor $F^{KL}$. In order to (re)-prove Theorem \ref{theorem_Lyubashenko}, we need to show that 
  \begin{equation}\label{eq_Lyubashenko}
   r(\phi^{-1}(x)) = F^{KL}((\phi,n)) r(x) F^{KL}((\phi, n))^{-1}, \quad \mbox{for all } x\in \mathcal{S}_{\zeta}(\mathbf{\Sigma}_{g,0}^*), (\phi, n) \in \widetilde{\Mod}(\Sigma_{g,1}).
   \end{equation}
  
  Let $\iota: \BT \to \mathcal{M}_c^{(1)}$ be the inclusion.  It is proved in \cite{KojuMurakami_QCharVar} that the following is a left Kan extension
  $$ \begin{tikzcd}
  \mathcal{M}_c^{(1)} \ar[r, "\mathcal{S}_{\zeta}^{[0]}"] & \mathcal{C} \\
  \BT \ar[u, "\iota"] \ar[ru, "\restriction{\mathcal{S}_{\zeta}^{[0]}}{\BT}"'] & {}
  \end{tikzcd}
  $$
  i.e. that $\mathcal{S}_{\zeta}^{[0]}$ is a (the) left Kan extension of its restriction to $\BT$ along $\iota$. Since the restriction of $F^{KL} \circ \widetilde{\partial} : \mathcal{M}_c^{(1)} \to \mathcal{C}$ to $\BT$ is equal to $\restriction{\mathcal{S}_{\zeta}^{[0]}}{\BT}$, by universal property of left Kan extension, we have a natural transformation $\eta:   \mathcal{S}_{\zeta}^{[0]} \Rightarrow F^{KL} \circ \widetilde{\partial}$. $\eta$ has the following down-to-earth description. Let $\mathcal{T}= [T,s] \in \mathcal{S}_{\zeta}^{[0]}(\mathbf{M})$ be the class of a stated tangle such that $T$ is a bottom tangle (such stated tangles span the whole $ \mathcal{S}_{\zeta}^{[0]}(\mathbf{M})$)  and suppose that $\partial \mathbf{M}=\Sigma_g^*$. 
  Let  $H_{g'}^* \subset \mathbf{M}$ be a tubular neighborhood of $T$ so that there exists $x_1, \ldots, x_{g'} \in \mathfrak{o}_q\SL_2$ such that the class $[T,s] \in \mathcal{S}_{\zeta}^{[0]}(H_{g'})$ equals $\mu^{(g'-1)}(x_1\otimes \ldots \otimes x_{g'})$. Consider the cobordism $C:= \mathbf{M}\setminus H_{g'}^* : \Sigma_{g'}^* \to \Sigma_g$. Then $\eta$ is characterized by the formula $\eta(\mathcal{T}) = F^{KL}(C, 0)\circ F^{KL}(\mu^{(g'-1)})(x_1\otimes \ldots \otimes x_{g'})$.   Now let $\phi \in \widetilde{\Mod}(\Sigma_{g,1})$ and let $C(\phi)$ be the associated mapping cylinder (pinched along $\partial \Sigma_g^*$). Then $\phi^{-1}(T)$ has a neighborhood $\phi^{-1}(H_{g'}) \subset \mathbf{M}$ with associated cobordism $C'= \mathbf{M}\setminus \phi^{-1}(H_{g'})= C(\phi) \circ C \circ C(\phi^{-1})$ so 
  $$\eta(\phi^{-1}(\mathcal{T}))= F^{KL}(C') \circ (\mu^{(g'-1)}(x_1 \otimes \ldots \otimes x_{g'})) = F^{KL}(\phi,n) \eta(\mathcal{T}) F^{KL}(\phi^{-1}, n).$$
  This implies Equation \eqref{eq_Lyubashenko} and thus provides an alternative conceptual proof of Theorem \ref{theorem_Lyubashenko}.

  \vspace{2mm}
  \par 
  \textit{Acknowledgments.} These are the notes of a series of talks given at a workshop held in Montpellier and Toulouse. The author thanks the organizers, S.Baseilhac, P.Roche and F.Costantino, for this very interesting workshop. 
 He acknowledges support from the European Research Council (ERC DerSympApp) under the European Union’s Horizon 2020 research and innovation program (Grant Agreement No. 768679).

\bibliographystyle{amsalpha}
\bibliography{biblio}

\end{document}